\documentclass[10pt]{article}

\usepackage{amsmath,amsthm}

\newtheorem{thm}{Theorem}[section]
\newtheorem{cor}[thm]{Corollary}

\newtheorem{lem}[thm]{Lemma}

\theoremstyle{remark}

\theoremstyle{plain}

\newtheoremstyle{note}
  {3pt}
  {3pt}
  {}
  {}
  {\itshape}
  {:}
  {.5em}
  {}

\theoremstyle{note}

\newtheoremstyle{citing}
  {3pt}
  {3pt}
  {\itshape}
  {}
  {\bfseries}
  {.}
  {.5em}
  {\thmnote{#3}}

\theoremstyle{citing}

\newtheoremstyle{break}
  {9pt}
  {9pt}
  {\itshape}
  {}
  {\bfseries}
  {.}
  {\newline}
  {}

\theoremstyle{break}

\theoremstyle{exercise}

\swapnumbers
\theoremstyle{plain}

\let\lvert=|\let\rvert=|

\begin{document}

\begin{center} 
\bf{Rational Generalised Moonshine from Abelian Orbifoldings of the Moonshine
Module}
\end{center}

\begin{center} 
Rossen Ivanov$^{1,2}$ and Michael Tuite$^{1,3}$\\ 

$^{1}$Department of Mathematical Physics,\\
National University of Ireland, Galway, Ireland\\

$^{2}$Institute for Nuclear Research and Nuclear Energy,\\
72 Tzarigradsko shosse, 1784 Sofia, Bulgaria\\

$^{3}$Dublin Institute for Advanced Studies,\\
10 Burlington Road, Dublin 4, Ireland\\

email: rossen.ivanov@nuigalway.ie, michael.tuite@nuigalway.ie
\end{center}

\begin{abstract}
We consider orbifoldings of the Moonshine Module with respect to the abelian
group generated by a pair of commuting Monster group elements with one of
prime order $p=2,3,5,7$ and the other of order $pk$ for $k=1$ or $k$ prime.
We show that constraints arising from meromorphic orbifold conformal field
theory allow us to demonstrate that each orbifold partition function with
rational coefficients is either constant or is a hauptmodul for an
explicitly found modular fixing group of genus zero. We thus confirm in the
cases considered the Generalised Moonshine conjectures for all rational
modular functions for the Monster centralisers related to the Baby Monster,
Fischer, Harada-Norton and Held sporadic simple groups. We also derive
non-trivial constraints on the possible Monster conjugacy classes to which
the elements of the orbifolding abelian group may belong.



Los Alamos Database Number: math.QA/0106027
\end{abstract}

\section{Introduction}

Orbifold constructions \cite{DHVW} in Meromorphic Conformal Field Theory 
\cite{Go},\cite{DGM} (MCFT) and Vertex Operator Algebras \cite{FLM2},\cite
{Ka},\cite{MN} provide the most natural setting for understanding Moonshine
phenomena \cite{T1},\cite{T2},\cite{DLM1}. The Moonshine Module \cite{FLM1},%
\cite{FLM2} whose automorphism group is the Monster finite sporadic group $%
{\bf M}$, is an orbifold MCFT constructed by orbifolding the Leech lattice
MCFT with respect to the group generated by a reflection involution. The
Moonshine Module partition function is the classical elliptic $J$ function
which is a hauptmodul for the genus zero modular group $SL(2,{\bf Z})$ and
is believed to be the unique MCFT with this partition function \cite{FLM2}.
Orbifolding the Moonshine Module with respect to the group generated by $%
g\in $ ${\bf M}$ leads naturally to the notion of the orbifold partition
function i.e. the Thompson series $T_{g}$. Monstrous Moonshine is mainly
concerned with the property, conjectured by Conway and Norton \cite{CN} and
subsequently proved by Borcherds \cite{B}, that each Thompson series $T_{g}$
is a hauptmodul for some genus zero fixing modular group. Assuming the
uniqueness of the Moonshine Module, this genus zero property is believed to
be equivalent to the following statement \cite{T2}: the only orbifold MCFT
that can arise by orbifolding with respect to $g$ is either the Moonshine
Module itself (for $g$ belonging to a so-called Fricke Monster conjugacy
class) or the Leech lattice MCFT (for $g$ belonging to a non-Fricke Monster
conjugacy class).

The Generalised Moonshine conjecture of Norton \cite{N1} is concerned with
modular functions associated with a commuting pair $g,h\in $ ${\bf M}$ and
asserts that each such modular function is either constant or is a
hauptmodul for some genus zero fixing group. No\ extension of the Borcherds'
approach to Monstrous Moonshine has yet been shown to be possible for
Generalised Moonshine. We argue that the most natural setting for these
conjectures is to consider orbifoldings of the Moonshine Module with respect
to the abelian group $\langle g,h\rangle $ generated by $g,h$ \cite{T3}. In
the cases where $\langle g,h\rangle $ can be generated by a single Monster
element, for example when $g$ and $h$ have coprime orders, then Generalised
Moonshine follows directly from Monstrous Moonshine. In this paper we
consider the case where $g$ is of prime order $p=2,3,5$ and $7$ and is of
Fricke type and $h$ is of order $pk$ for $k=1$ or $k$ prime. We confirm
Norton's conjecture for modular functions with rational coefficients in
these cases by considering orbifold modular properties and some consistency
conditions arising from the orbifolding procedure. We also demonstrate a
number of other non-trivial aspects of Generalised Moonshine such as
properties of the character expansion of Generalised Moonshine functions in
terms of irreducible characters for the centraliser of $g$ in ${\bf M}$ and
constraints on the possible Monster conjugacy classes to which the elements
of $\langle g,h\rangle \,$ may belong.

We begin in Section 2 with a general review of Abelian orbifold
constructions in Meromorphic Conformal Field Theory (MCFT). We also briefly
review the construction of the Moonshine Module and the relationship between
the genus zero property of Thompson series $T_{g}$ in Monstrous Moonshine
for $g\in $ ${\bf M}$ and evidence to support the claim that the only
possible MCFTs obtainable by orbifolding the Moonshine Module with respect
to the group generated by $g$ are the Moonshine Module, for $g$ Fricke, and
the Leech lattice MCFT, for $g$ non-Fricke.

In Section 3 we begin with a discussion of general properties for
Generalised Moonshine Functions (GMF) following from the orbifold
considerations of Section 2. We then prove two theorems concerning
constraints that arise from the consistency of orbifolding the Moonshine
Module with respect to $\langle g,h\rangle $ under various choices of
generators. These constraints are exploited in Section 4 in order to
determine the residues of singular cusps of GMFs. In particular if all
elements of $\langle g,h\rangle $ are non-Fricke the GMF is constant. We
then prove some modular properties of GMFs with rational coefficents in the
cases where $g$ is Fricke of prime order $p=2,3,5$ and $7$ and $h$ is of
order $pk$ for $k=1$ or $k$ prime. This analysis in part relies on
properties of the characters of the centralisers of the Monster related to
the Baby Monster, Fischer, Harada-Norton and Held sporadic simple groups. We
also highlight the importance of the Monster conjugacy classes to which the
elements of $\langle g,h\rangle $ belong in determining the possible
singularities of a GMF.

In Section 4 we give a comprehensive analysis of the possible singularity
structure of GMFs for the cases under consideration. In each case we
demonstrate that either the given singularity structure is inconsistent or
else all singularities of the GMF can be identified under some genus zero
fixing group for which the GMF is a hauptmodul. Thus we obtain non-trivial
constraints on the Monster classes to which the elements of $\langle
g,h\rangle $ may belong and verify the Generalised Moonshine conjecture in
these cases.

In Appendix A we review the definitions of standard modular groups. In
Appendix B we consider the first ten coefficients of a general GMF as
character expansions for the Monster centralisers related to the Baby
Monster and Fischer centraliser subgroups.

\section{Abelian Orbifolds and Monstrous Moonshine}

\subsection{Self-Dual $C=24$ Meromorphic CFTs}

A Meromorphic CFT (MCFT) or chiral CFT is a CFT whose n-point functions are
all meromorphic as described in \cite{Go}, \cite{DGM}. A\ MCFT essentially
corresponds to a Vertex Operator Algebra in the pure mathematics literature
as reviewed in \cite{FLM2}, \cite{Ka}, \cite{MN}. We briefly review some of
their basic properties. A\ MCFT $(\mathcal{V},\mathcal{H)}$ consists of a
Hilbert space of states $\mathcal{H\,}$together with a set of vertex
operators $\mathcal{V}\equiv \{V(\psi ,z)=\sum_{n\in {\bf Z}}\psi
_{n}z^{-n-h_{\psi }}|\psi \in \mathcal{H}\}$ where $h_{\psi }$ is the
conformal weight (see below) and where each mode $\psi _{n}\,$acts as a
linear operator on $\mathcal{H}$. The vertex operators are local meaning
that for given $\varphi ,\psi \in \mathcal{H}$ and sufficiently large $N$ 
\begin{equation}
(V(\psi ,z)V(\varphi ,w)-V(\varphi ,w)V(\psi ,z))(z-w)^{N}=0.  \label{local}
\end{equation}
$\mathcal{H}$ contains a distinguished vacuum state $|0\rangle \,$such that 
\begin{eqnarray}
V(|0\rangle ,z) &=&\mathrm{Id}_{\mathcal{H}},  \label{identity} \\
{\mathrm{lim}}_{z\rightarrow 0}V(\psi ,z)|0\rangle &=&\psi .
\label{creation}
\end{eqnarray}
$\mathcal{H}$ also contains a Virasoro state $\omega $ with vertex operator
Laurant expansion $V\left( \omega ,z\right) =\sum_{n\in {\bf Z}%
}L_{n}z^{-n-2}\,$where $L_{n}\,$ generates the Virasoro algebra of central
charge $C$ 
\begin{equation}
\lbrack L_{m},L_{n}]=(m-n)L_{m+n}+\frac{1}{12}C(m^{3}-m)\delta _{m,-n}.
\label{Virasoro}
\end{equation}
We will consider here MCFTs of central charge $C=24$ only. The vector space $%
\mathcal{H}$ is decomposed into finite dimensional spaces $\mathcal{H}_{n}$
with non-negative integral $L_{0}$ grading $n$, the conformal weight. Taking
these various properties together, the vertex operators then satisfy the
Operator Product Expansion (OPE) for $|z|>|w|$%
\begin{eqnarray}
V(\psi ,z)V(\varphi ,w) &=&V(V(\psi ,z-w)\varphi ,w)  \label{OPE} \\
&=&{\sum}_{n\geq 0}{\sum }_{\chi \in \mathcal{H}_{n}}%
C_{\psi \varphi }^{\chi }V(\chi ,w)(z-w)^{n-h_{\psi }-h_{\varphi }}
\label{OPE physics}
\end{eqnarray}
with $\psi _{h_{\varphi }-n}(\varphi )=\sum_{\chi \in \mathcal{H}%
_{n}}C_{\psi \varphi }^{\chi }\chi $ and where the sum is taken over some
basis for $\mathcal{H}_{n}$ \cite{Go}, \cite{DGM}.

The automorphism group $\mathrm{Aut}(\mathcal{V})$ of $\mathcal{V}$ is the
group of linear transformations $g:\mathcal{H\rightarrow H}$ which preserves
the Virasoro state $\omega $ and where 
\begin{equation}
gV(\psi ,z)g^{-1}=V(g\psi ,z).  \label{Auto of MCFT}
\end{equation}
The OPE (\ref{OPE}) is then invariant under $\mathrm{Aut}(\mathcal{V})$. In
the case of the Moonshine Module $\mathcal{V}^{\natural }$, $\mathrm{Aut}(%
\mathcal{V}^{\natural })={\bf M}$ the Monster group, which has a unitary
action on $\mathcal{V}^{\natural }$. We will assume that $\mathrm{Aut}(%
\mathcal{V})$ has a unitary action from now on (which is expected physically
from unitarity).

A\ representation $(\mathcal{U},\mathcal{K)}$ of a MCFT $(\mathcal{V},%
\mathcal{H)}$ consists of a vector space (module) $\mathcal{K\ }$together
with a set of local vertex operators $\mathcal{U\equiv }\{U(\psi ,z)|\psi
\in \mathcal{H}\}$ where the modes of $U(\psi ,z)\,$act as linear operators
on $\mathcal{K}$ with $U(|0\rangle ,z)=\mathrm{Id}_{\mathcal{K}}.$ These
operators satisfy the OPE 
\begin{equation}
U(\psi ,z)U(\varphi ,w)=U(V(\psi ,z-w)\varphi ,w),\quad |z|>|w|.
\label{rep of MCFT}
\end{equation}
$(\mathcal{U},\mathcal{K})\mathcal{\,}$is an irreducible representation if $%
\mathcal{K}$ contains no non-trivial submodule invariant under the modes of $%
\{U(\psi ,z)\}$. $(\mathcal{V},\mathcal{H)}$ is said to be a Self-Dual MCFT
or a Holomorphic VOA if $(\mathcal{V},\mathcal{H)}$ is the unique
irreducible representation for itself. We assume that $\mathcal{K\,}$ is
decomposed into Verma modules of the Virasoro algebra of $\mathcal{V\,}$
with non-negative integral $L_{0}$ grading i.e. we consider unitary Virasoro
representations. Let $\mathcal{K}_{0}\subset \mathcal{K\,}$ denote the
subspace of lowest $L_{0}$ grading. Then for an irreducible representation, $%
\mathcal{K}$ is generated by the action of the modes of $\{U(\psi ,z)\}$ on $%
\mathcal{K}_{0}$.

The automorphism group $\mathrm{Aut}(\mathcal{U})$ of $\mathcal{U}$ is the
group of linear transformations $\hat{g}:\mathcal{K\rightarrow K}$ of the
form 
\begin{equation}
\hat{g}U(\psi ,z)\hat{g}^{-1}=U(g\psi ,z)  \label{Auto of rep}
\end{equation}
$\,$for some $g\in \mathrm{Aut}(\mathcal{V})$. A general element of $\mathrm{%
Aut}(\mathcal{V})\,$ may give rise to a mapping between different
representations and so we denote by $\mathrm{Aut}_{\mathcal{K}}(\mathcal{V})$
the subgroup of $\mathrm{Aut}(\mathcal{V})$ associated with (\ref{Auto of
rep}). Clearly $\mathrm{Aut}_{\mathcal{K}}(\mathcal{V})$ acts projectively
on $\mathcal{K}$ with a natural homomorphism from $\mathrm{Aut}(\mathcal{U})$
to $\mathrm{Aut}_{\mathcal{K}}(\mathcal{V})$ whose kernel consists of phase
multipliers assuming that $\mathrm{Aut}(\mathcal{U})\,$ acts unitarily on $%
\mathcal{K}$ i.e. $\mathrm{Aut}(\mathcal{U})=U(1).\mathrm{Aut}_{\mathcal{K}}(%
\mathcal{V})$. For an irreducible representation, where $\mathcal{K}_{0}$ is
one dimensional, then the $U(1)\,$ subgroup of $\mathrm{Aut}(\mathcal{U})$
has the same action on all elements of $\mathcal{K}\,$ so that 
\begin{equation}
\mathrm{Aut}(\mathcal{U})=U(1)\times \mathrm{Aut}_{\mathcal{K}}(\mathcal{V}%
),\quad  \label{aut(V) for dim(k0)=1}
\end{equation}
if\textrm{\ }$\dim (\mathcal{K}_{0})=1$.

Let us now assume that $\mathcal{V\,}$ is Self-Dual. The characteristic
function (or genus one partition function) for $\mathcal{V}$ of central
charge $C=24$ is given by

\begin{equation}
Z(\tau )=\mathrm{Tr}_{\mathcal{H}}(q^{L_{0}-1}),\text{ }q=e^{2\pi i\tau }
\label{PF}
\end{equation}
where $\tau \in {\bf H}$, the upper half complex plane, is the usual
elliptic modular parameter. $\mathcal{V}$ has integral grading so that $%
Z(\tau )$ $\,$is invariant under $T:\tau \rightarrow \tau +1$. Self-duality
implies that $Z(\tau )$ $\,$is invariant under $S:\tau \rightarrow -1/\tau $
so that $Z(\tau )$ is invariant under the modular group {$SL$}$(2,{\bf Z})$,
generated by $S,T$ \cite{Z}. Hence $Z(\tau )$ is uniquely determined up to
an additive constant by $J(\tau )$, the hauptmodul for {$SL$}$(2,{\bf Z})$ 
\cite{Se} 
\begin{eqnarray}
Z(\tau ) &=&J(\tau )+N_{0}\text{, }  \nonumber \\
J(\tau ) &=&\frac{E_{4}^{3}}{{\eta ^{24}}}-744{=}\frac{1}{q}+0+196884q+...,
\label{J}
\end{eqnarray}
where $\eta (\tau )=q^{1/24}\prod_{n>0}(1-q^{n})$ is the Dedekind eta
function, $E_{n}(\tau )$ is the Eisenstein form of weight $n$ \cite{Se} and $%
N_{0}$ is the number of conformal weight 1 operators in $\mathcal{V}$.
Examples of such theories are lattice models, which we denote by $\mathcal{V}%
^{\Lambda }$, where $\Lambda $ is one of the Niemeier even self-dual 24
dimensional lattices. Then $Z(\tau )=\Theta _{\Lambda }/\eta ^{24}$ with $%
\Theta _{\Lambda }=\sum_{\lambda \in \Lambda }q^{\lambda ^{2}/2}$ the
lattice theta function for $\Lambda $. In this paper we will be particularly
concerned with the Leech lattice for which $N_{0}=24$ and the Moonshine
Module $\mathcal{V}^{\natural }$ for which $N_{0}=0$.

\subsection{Twisted Sectors of a Self-Dual MCFT}

In this section we briefly review relevant aspects of twisted representation
of a MCFT. Let $G=\mathrm{Aut}(\mathcal{V})$, denote the automorphism group (%
\ref{Auto of MCFT}) for $\mathcal{V}$. Consider $g\in G$ of finite order $%
o(g)=$ $n$ and define the trace function 
\begin{equation}
Z\left[ 
\begin{array}{c}
g \\ 
1
\end{array}
\right] (\tau )\equiv \mathrm{Tr}_{\mathcal{H}}(gq^{L_{0}-1}).  \label{Zg1}
\end{equation}
Clearly $Z\left[ 
\begin{array}{c}
1 \\ 
1
\end{array}
\right] (\tau )\equiv Z(\tau )$ of (\ref{PF}). In the case where $\mathcal{V}%
=\mathcal{V}^{\natural }$, the Moonshine Module, (\ref{Zg1}) is the Thompson
series (see below). Let $\mathcal{H}^{(j)}$ denote the eigenspace of $%
\mathcal{H}$ for $g$ with eigenvalue $\omega _{n}^{j}\,$for $\omega
_{n}\equiv \exp (2\pi i/n)$. The twisted representation $(\mathcal{V}_{g},%
\mathcal{H}_{g})$ of $(\mathcal{V},\mathcal{H})\,\,$consists of a set of
vertex operators $\mathcal{V}_{g}$ $\equiv \{V_{g}(\psi ,z)|\psi \in 
\mathcal{H}\}$ with mode expansion 
\begin{equation}
V_{g}(\psi ,z)={\sum }_{m\in {\bf Z}+j/n}\tilde{\psi}%
_{m}z^{-m-1},\quad \psi \in \mathcal{H}^{(j)}  \label{Vgmodes}
\end{equation}
whose modes $\tilde{\psi}_{m}\,$ are linear operators on $\mathcal{H}_{g}\,$%
. The twisting property corresponds to the monodromy relation $V_{g}(\psi
,e^{2\pi i}z)=V_{g}(g^{-1}\psi ,z).$ Furthermore the twisted vertex
operators satisfy the `twisted` non-meromorphic OPE 
\begin{equation}
V_{g}(\psi ,z)V_{g}(\varphi ,w)=\omega _{n}^{-jb}V_{g}(V(\psi ,z-w)\varphi
,w),\quad  \label{twisted OPE}
\end{equation}
for $\psi \in \mathcal{H}^{(j)}$ and where $b=0,1,2...,n-1\,$labels the
sheet for the branched $n$-fold covering for $(w/z)^{1/n}$. For $(\mathcal{V}%
,\mathcal{H})$ self-dual, $(\mathcal{V}_{g},\mathcal{H}_{g})\,$ always
exists and is unique up to isomorphism \cite{DLM1}. Furthermore, under
conjugation by any element $x\in G$ then $x(\mathcal{V}_{g})x^{-1}$ is
isomorphic to $\mathcal{V}_{xgx^{-1}}.$

Considering (\ref{twisted OPE}) for $\psi ,\varphi \in \mathcal{H}^{(0)}$ we
see that $(\mathcal{V}_{g},\mathcal{H}_{g})$ forms a reducible
representation for $(\mathcal{V}^{(0)},\mathcal{H}^{(0)})\,\,$as in (\ref
{rep of MCFT}) where $\mathcal{V}^{(0)}\,$are the vertex operators for $%
\mathcal{H}^{(0)}$. $\mathcal{H}_{g}\ $can therefore be decomposed into $%
L_{0}$ eigenspaces $\mathcal{H}_{g}=\bigoplus_{m=0}^{\infty }\mathcal{H}%
_{g,m}$ where $\mathcal{H}_{g,m}$ has rational $L_{0}$ eigenvalue $%
E_{0}^{g}+1+m/n$. The subspace $\{\sigma _{g}^{a}\}$ for $\,a=1,2,...,D_{g}$
of lowest grade $E_{0}^{g}+1\ $is called the $g$-twisted vacuum space. $%
E_{0}^{g}$ is called the vacuum energy and $D_{g}$ is called the vacuum
degeneracy.

The automorphism group, $\mathrm{Aut}(\mathcal{V}_{g})$, preserving (\ref
{twisted OPE}) can be defined in a way similar to (\ref{Auto of rep}) as an
extension of $C_{g}=\{h\in G|gh=hg\}$, the centraliser of $g$ in $G$. Since $%
\mathcal{V}_{g}$ $\,$is unique for a self-dual MCFT we have $\mathrm{Aut}(%
\mathcal{V}_{g})=U(1).C_{g}$. If the twisted vacuum is unique ($D_{g}=1$)
then $\mathrm{Aut}(\mathcal{V}_{g})=U(1)\times C_{g}$ from (\ref{aut(V) for
dim(k0)=1}).

Let $\hat{g}\in \mathrm{Aut}(\mathcal{V}_{g})$ denote the lifting of $g$
with action on the twisted vacuum as follows 
\begin{equation}
\hat{g}\sigma _{g}^{a}=\exp (-2\pi iE_{0}^{g})\sigma _{g}^{a},
\label{ghat on twisted vacuum}
\end{equation}
where $\exp (-2\pi iE_{0}^{g})\,$is a $U(1)$ phase. Then (\ref{Vgmodes})
implies that in general 
\begin{equation}
\hat{g}\psi _{g}=\exp (-2\pi ih_{g})\psi _{g},
\label{ghat on all twisted states}
\end{equation}
where $\psi _{g}\in $ $\mathcal{H}_{g}$ has Virasoro grading $h_{g}$.
Clearly $n|o(\hat{g})$, where $o(\hat{g})$ is the order of $\hat{g}$, since $%
\hat{g}^{n}\,$is a lifting of the identity element of $C_{g}$ so that $%
E_{0}^{g}\in {\bf Z}/o(\hat{g})$. We say that $g$ is a Normal element of $G$
if $nE_{0}^{g}\in {\bf Z}$ so that $\hat{g}$ is of order $n$, otherwise we
say that $g$ is an Anomalous element of $G$.

Let $g$ be a normal element of $G$ and consider the twisted spaces $(%
\mathcal{V}_{g^{k}},\mathcal{H}_{g^{k}})$ for $k=1,2,...,n-1$. Let $\mathcal{%
H}_{g^{k}}^{(j)}$ denote the eigenspace of $\mathcal{H}_{g^{k}}\,$with $\hat{%
g}$ eigenvalue $\omega _{n}^{j}$ where $\mathcal{H}_{g^{k}}^{(j)}$ can be
further decomposed into $L_{0}$ eigenspaces $\mathcal{H}_{g^{k},m}^{(j)}$.
Then $\{(\mathcal{V}_{g^{k}}^{(j)},\mathcal{H}_{g^{k}}^{(j)})\}$ comprises
the $n^{2}$ irreducible representations for the MCFT $(\mathcal{V}^{(0)},%
\mathcal{H}^{(0)})$ \cite{DLM1}. Clearly $\mathrm{Aut}(\mathcal{V}%
^{(0)})\supseteq G_{g}$ where $G_{g}=C_{g}/\langle g\rangle $ so that $%
\mathrm{Aut}(\mathcal{V}_{g}^{(j)})=U(1).G_{g}$. Note that $\mathrm{Aut}(%
\mathcal{V}_{g}^{(j)})$ depends on $j$ in general. In particular, if the
twisted vacuum is unique then $\sigma _{g}^{1}\in \mathcal{H}_{g}^{(j_{0})}$
for some $j_{0}$ and $\mathrm{Aut}(\mathcal{V}_{g}^{(j_{0})})=U(1)\times
G_{g}$ from (\ref{aut(V) for dim(k0)=1}) where $\omega _{n}^{j_{0}}=$ $\exp
(-2\pi iE_{0}^{g})$.

We next define the trace function for $\mathcal{H}_{g}=\bigoplus_{j=0}^{n-1}$
$\bigoplus_{m=0}^{\infty }\mathcal{H}_{g,m}^{(j)}$ as follows 
\begin{eqnarray}
Z\left[ 
\begin{array}{c}
1 \\ 
g
\end{array}
\right] (\tau ) &\equiv &\mathrm{Tr}_{\mathcal{H}_{g}}(q^{L_{0}-1})
\label{Z(1,g)} \\
&=&\sum_{j=0}^{n-1}\mathrm{Tr}_{\mathcal{H}%
_{g}^{(j)}}(q^{L_{0}-1})=D_{g}q^{E_{0}^{g}}+..., \\
\mathrm{Tr}_{\mathcal{H}_{g}^{(j)}}(q^{L_{0}-1}) &=&q^{-j/n}
\sum_{m=0}^\infty D_{g,m}^{(j)}q^{m}.  \label{Trace H^(j)_g}
\end{eqnarray}
where the coefficient $D_{g,m}^{(j)}$ is the dimension of the representation 
$\rho _{g,m}^{(j)}$ of $\mathrm{Aut}(\mathcal{V}_{g}^{(j)})$ defined by $%
\mathcal{H}_{g,m}^{(j)}$. Then $D_{g}\,$is the dimension of $\rho _{g}^{0}$,
the representation of $\mathrm{Aut}(\mathcal{V}_{g}^{(j)})$ acting on the
twisted vacuum.

We can similarly define the general trace function for $\hat{h}\in \mathrm{%
Aut}(\mathcal{V}_{g})$ lifted from $h\in C_{g}$ by

\begin{eqnarray}
Z\left[ 
\begin{array}{c}
h \\ 
g
\end{array}
\right] (\tau ) &\equiv &\mathrm{Tr}_{\mathcal{H}_{g}}(\hat{h}q^{L_{0}-1})
\label{Z(h,g)} \\
&=&\sum_{j=0}^{n-1}q^{-j/n}
\sum_{m=0}^\infty\chi _{g,m}^{(j)}(\hat{h})q^{m}=\chi _{g}^{0}(\hat{h}%
)q^{E_{0}^{g}}+...,  \label{Z(h,g) char expansion}
\end{eqnarray}
where $\chi _{g,m}^{(j)}(\hat{h})=\mathrm{Tr}(\rho _{g,m}^{(j)}(\hat{h}%
))\,\, $denotes a character of $\hat{h}\in \mathrm{Aut}(\mathcal{V}%
_{g}^{(j)})$ and $\chi _{g}^{0}$ is the character for $\rho _{g}^{0}$. We
assume below that a particular choice for $\hat{h}\,$ can be made which
resolves the ambiguity inherent in the notation $Z\left[ 
\begin{array}{c}
h \\ 
g
\end{array}
\right] \,$ denoting the trace (\ref{Z(h,g)}). When $D_{g}=1$, then $\,\chi
_{g,m}^{(j)}(\hat{h})/\chi _{g}^{0}(\hat{h})$ is a character for $C_{g}$ for 
$j\neq j_{0}$ and is a character for $G_{g}\,$ for $j=j_{0}$ where $\omega
_{n}^{j_{0}}=$ $\exp (-2\pi iE_{0}^{g})$.

For general commuting elements $g,h$, the trace function (\ref{Z(h,g)})
transforms under a modular transformation with respect to $\gamma =\left( 
\begin{array}{cc}
a & b \\ 
c & d
\end{array}
\right) \in SL(2,{\bf Z)}$ in the following way

\begin{equation}
Z\left[ 
\begin{array}{c}
h \\ 
g
\end{array}
\right] (\tau )=\varepsilon (g,h;\gamma )Z\left[ 
\begin{array}{c}
h \\ 
g
\end{array}
\right] ^{\gamma }(\gamma \tau ),  \label{SLT with phase}
\end{equation}
where $\gamma \tau =\frac{a\tau +b}{c\tau +d}$, $\left[ 
\begin{array}{c}
h \\ 
g
\end{array}
\right] ^{\gamma }\equiv \left[ 
\begin{array}{c}
h^{a}g^{b} \\ 
h^{c}g^{d}
\end{array}
\right] $ and where $\varepsilon (g,h;\gamma )$ is a phase multiplier \cite
{DHVW}, \cite{DLM1}. We will assume that if all elements in $\langle
g,h\rangle $, the group generated by $g$ and $h$, are normal then we may
choose a lifting of $h^{a}g^{b}\,$ such that the phase multiplier $%
\varepsilon (g,h;\gamma )=1\,$ i.e. there are no global phase anomalies \cite
{Va}. We will employ the abbreviation $h$ for the chosen lifting $\hat{h}\in 
\mathrm{Aut}(\mathcal{V}_{g})$ of $h\in C_{g}$ from now on. In general, the
absence of this phase multiplier is an essential ingredient in the
orbifolding procedure that we discuss below.

Suppose that the $g$ twisted vacuum is one dimensional i.e. $D_{g}=1$. Let $%
\phi _{g}(h)\,$ denote the lifting of $h$ in its action on this twisted
vacuum (and hence giving the extension of $h$ on all of $\mathcal{H}_{g}$)
where in particular $\phi _{g}(1)=1$ and $\phi _{g}(g)=\exp (-2\pi
iE_{0}^{g})$. We conjecture that $\phi _{g}(h)\,\,$has the following
properties where $g^{a}h^{b}$ is a normal class for all $a,b$: 
\begin{eqnarray}
\phi _{g}(g^{a}h^{b}) &=&\phi _{g}(g)^{a}\phi _{g}(h)^{b},
\label{phi_gh_gahb} \\
\phi _{g}(h) &\in &\langle \omega _{n}\rangle .\text{ }  \label{phi_gh_n}
\end{eqnarray}
(\ref{phi_gh_gahb}) follows from modular invariance with $Z\left[ 
\begin{array}{c}
h^{b} \\ 
g
\end{array}
\right] (\tau +1)=Z\left[ 
\begin{array}{c}
g^{-1}h^{b} \\ 
g
\end{array}
\right] (\tau )$ $\,$using (\ref{ghat on all twisted states}) so that $\phi
_{g}(g^{a}h^{b})=\phi _{g}(g)^{a}\phi _{g}(h^{b})$ and the assumption that $%
\widehat{g^{a}h^{b}}=(\hat{g})^{a}(\hat{h})^{b}$ for normal classes $%
g^{a}h^{b}$. We will prove (\ref{phi_gh_n}) assuming (\ref{phi_gh_gahb}) for
the specific examples of Generalised Moonshine Functions that we consider
later on.

\subsection{Orbifolding a MCFT}

Assume that all elements of $\langle g\rangle \simeq {\bf Z}_{n}$, the
abelian group of order $n$ generated by $g$, are normal elements of $\mathrm{%
Aut}(\mathcal{V})$. Then $(\mathcal{V}^{(0)},\mathcal{H}^{(0)})$ has $n^{2}$
irreducible representations $(\mathcal{V}_{g^{k}}^{(j)},\mathcal{H}%
_{g^{k}}^{(j)})$ for $j,k=0...,n-1$. The $\langle g\rangle $ orbifold MCFT $%
\mathcal{V}_{\text{orb}}^{\langle g\rangle }$ is the MCFT with Hilbert space 
$\mathcal{H}_{\text{orb}}^{\langle g\rangle }\equiv \oplus _{k=0}^{n-1}%
\mathcal{H}_{g^{k}}^{(0)}$ . We assume that we can augment the operators $%
\mathcal{V}^{(0)}$ with appropriate local operators so that the OPE (\ref
{OPE}) is satisfied. The characteristic function of $\mathcal{V}_{\text{orb}%
}^{\langle g\rangle }\,$is $Z_{\text{orb}}^{\langle g\rangle }=\frac{1}{n}%
\sum_{l,k=0}^{n-1}Z\left[ 
\begin{array}{c}
g^{l} \\ 
g^{k}
\end{array}
\right] $ which is a modular invariant from (\ref{SLT with phase}) since $%
\varepsilon (g,h;\gamma )\equiv 1$. Hence $\mathcal{V}_{\text{orb}}^{\langle
g\rangle }$ is a self-dual MCFT and so $Z_{\text{orb}}^{\langle g\rangle
}(\tau )=J(\tau )+N_{0}^{\langle g\rangle }$ as in (\ref{J}) where $%
N_{0}^{\langle g\rangle }$ is the number of conformal weight 1 operators in $%
\mathcal{V}_{\text{orb}}^{\langle g\rangle }$.

We can similarly consider the orbifold MCFT $\mathcal{V}_{\text{orb}%
}^{\langle g,h\rangle }\,$ found by orbifolding with respect to the abelian
group $\langle g,h\rangle $ of order $|\langle g,h\rangle |$ generated by
two commuting elements $g,h$ where all the elements of $\langle g,h\rangle $
are assumed to be normal. $\mathcal{V}_{\text{orb}}^{\langle g,h\rangle }$
has Hilbert space $\mathcal{H}_{\text{orb}}^{\langle g,h\rangle }=\mathcal{P}%
_{\langle g,h\rangle }(\oplus _{v\in \langle g,h\rangle }\mathcal{H}_{v})$
where $\mathcal{P}_{\langle g,h\rangle }\equiv \frac{1}{|\langle g,h\rangle |%
}\sum_{v\in \langle g,h\rangle }v$ denotes the projection with
respect to the group $\langle g,h\rangle $. Again we assume that we may
augment the MCFT $\mathcal{P}_{\langle g,h\rangle }\mathcal{V}$ with
appropriate local vertex operators. $\mathcal{V}_{\text{orb}}^{\langle
g,h\rangle }$ then has modular invariant characteristic function $Z_{\text{%
orb}}^{\langle g,h\rangle }=\frac{1}{|\langle g,h\rangle |}\sum_{u,v\in
\langle g,h\rangle }Z\left[ 
\begin{array}{c}
u \\ 
v
\end{array}
\right] $. $\mathcal{V}_{\text{orb}}^{\langle g,h\rangle }$ is therefore a
self-dual MCFT and so $Z_{\text{orb}}^{\langle g,h\rangle }(\tau )=J(\tau
)+N_{0}^{\langle g,h\rangle }$.

We assume that these orbifold MCFTs further can be considered as various
embeddings in a larger non-meromorphic CFT $(\mathcal{V}^{\prime },\mathcal{H%
}^{\prime })$ with Hilbert space $\mathcal{H}^{\prime }=\oplus _{v\in
\langle g,h\rangle }\mathcal{H}_{v}$. In this CFT, all twisted states are
created by vertex operators satisfying some non-meromorphic OPE of the
generic form

\begin{equation}
V(\psi ,z)V(\varphi ,w)=\sum_\chi C_{\psi \varphi }^{\chi
}V(\chi ,w)(z-w)^{h_{\chi }-h_{\psi }-h_{\varphi }},
\label{Nonmermorphic OPE}
\end{equation}
for $\psi \in \mathcal{H}_{g},\quad \varphi \in \mathcal{H}_{h}$ and $
\chi \in \mathcal{H}_{gh}$ and similarly for all commuting pairs in $\langle g,h\rangle $. $\mathcal{V}%
_{\text{orb}}^{\langle g,h\rangle }\,$then consists of all $\mathcal{V}%
^{\prime }$ operators invariant under $\langle g,h\rangle $. A rigorous
discussion in the case of the reflection automorphism orbifolding of the
Leech lattice theory appears in \cite{Hu}.

Consider independent generators $g,h$ of $\langle g,h\rangle \,$i.e. $%
g^{A}\neq h^{B}$ for all $A=1,...,o(g)-1$ and $B=1,...,o(h)-1$ where $o(g)$
is the order of $g$ etc. Then $|\langle g,h\rangle |=o(g)o(h)$ and $\mathcal{%
P}_{\langle g,h\rangle }=\mathcal{P}_{g}\mathcal{P}_{h}$, where $\mathcal{P}%
_{g}=\frac{1}{o(g)}\sum_{k=0}^{o(g)-1}g^{k}$ denotes the projection operator
with respect to $g$. The orbifold $\mathcal{V}_{\text{orb}}^{\langle
g,h\rangle }\,$can then be considered as a composition of orbifoldings for
any independent generators $g,h$ where 
\begin{eqnarray}
\mathcal{V}_{\text{orb}}^{\langle g,h\rangle } &=&\mathcal{P}%
_{g}(\sum_{k=0}^{o(g)-1}\mathcal{P}_{h}(\sum_{l=0}^{o(h)-1}\mathcal{V}%
_{g^{k}h^{l}}))  \nonumber \\
&=&\mathcal{P}_{g}(\sum_{k=0}^{o(g)-1}(\mathcal{V}_{\text{orb}}^{\langle
h\rangle })_{g^{k}}))=(\mathcal{V}_{\text{orb}}^{\langle h\rangle })_{\text{%
orb}}^{\langle g\rangle },  \label{Hghorb}
\end{eqnarray}
using the uniqueness of the twisted sectors for a self-dual MCFT and the
assumed embedding of $\mathcal{V}_{\text{orb}}^{\langle g,h\rangle }$, $%
\mathcal{V}_{\text{orb}}^{\langle g\rangle }$ and $\mathcal{V}_{\text{orb}%
}^{\langle h\rangle }$ in $\mathcal{V}^{\prime }$. Later on we consider the
consistency of (\ref{Hghorb}) under the various possible independent choices
for the generators of $\langle g,h\rangle $ in proving Theorems \ref
{theorem_proj1} and \ref{theorem_proj2} in Section 3.

\subsection{The Moonshine Module and Monstrous Moonshine}

The Moonshine module $\mathcal{V}^{\natural }$ is historically the first
example of a self-dual orbifold MCFT \cite{FLM1} and is constructed as a $%
{\bf Z}_{2}$ orbifolding of $\mathcal{V}^{\Lambda }$, which denotes the
Leech lattice MCFT from now on. The ${\bf Z}_{2}$ automorphism $r\in \mathrm{%
Aut}(\mathcal{V}^{\Lambda })$ is a lifting of the lattice $\,$reflection
symmetry chosen so that $\mathcal{P}_{r}\mathcal{H}^{\Lambda }$ contains no
Virasoro level one states. The $r$-twisted space $\mathcal{H}_{r}^{\Lambda }$
has vacuum energy $E_{0}^{r}=1/2>0$ ($r$ is a normal element of $\mathrm{Aut}%
(\mathcal{V}^{\Lambda })$) and hence contains no Virasoro level one states.
The resulting orbifold MCFT, $\mathcal{V}^{\natural }\equiv (\mathcal{V}%
^{\Lambda })_{\mathrm{orb}}^{\langle r\rangle }$, therefore has
characteristic function $J(\tau )$ (\ref{J}) with $\mathrm{Aut}(\mathcal{V}%
^{\natural })={\bf M}$, the Monster group \cite{FLM1,FLM2}. We can identify
a 'dual' automorphism $r^{*}\in {\bf M}$ where $\mathcal{P}_{r}(\mathcal{H}%
^{\Lambda })$ (respectively $\mathcal{P}_{r}(\mathcal{H}_{r})$) is even
(respectively odd) under $r^{*}$ \cite{T2}. This is an obvious automorphism
of the non-meromorphic OPE (\ref{Nonmermorphic OPE}) for $\mathcal{V}%
^{\prime }$ for $g,h\in \langle r\rangle $. Then orbifolding $\mathcal{V}%
^{\natural }$ with respect to $r^{*}$ we recover $\mathcal{V}^{\Lambda }$.
Furthermore, one obtains the Monster centraliser $C_{r^{*}}=2_{+}^{1+24}.%
\mathrm{Co}_{1}$ where $\mathrm{Co}_{1}$ denotes the Conway simple group and 
$2_{+}^{1+24}$ is an extra-special 2-group \cite{Gr,FLM1,FLM2}.

It is conjectured that $\mathcal{V}^{\natural }$ is characterised (up to
isomorphism) as the unique self-dual $C=24$ MCFT with characteristic
function $J(\tau )$ \cite{FLM2}. We may consider other ${\bf Z}_{n}$
orbifoldings of $\mathcal{V}^{\Lambda }$ with characteristic function $%
J(\tau )$ which should reproduce $\mathcal{V}^{\natural }$ according to this
conjecture. In general, we can classify all automorphisms $a\in \mathrm{Aut}(%
\mathcal{V}^{\Lambda })$ lifted from automorphisms $\overline{a}\in \mathrm{%
Co}_{0}$ the Leech lattice automorphism group for which $\mathcal{V}%
_{a}^{\Lambda }$ can be explicitly constructed satisfying the following
constraints \cite{T2}

{(i)} $\mathcal{P}_{a}\mathcal{H}^{\Lambda }$ contains no Virasoro level one
states i.e. $\bar{a}$ $\,$is fixed point free.

{(ii) }$\mathcal{H}_{a}^{\Lambda }$ is non-tachyonic (i.e. $E_{0}^{a}\geq 0$
) and furthermore contains no Virasoro level one states so that 
\begin{equation}
E_{0}^{a}>0  \label{E0a}
\end{equation}

{(iii}) $a$ is a normal element of $\mathrm{Aut}(\mathcal{V}^{\Lambda })$.

There are 38 classes of $\mathrm{Co}_{0}$ obeying these constraints
including the 5 prime ordered cases considered by Dong and Mason \cite{DM}.

For each of these 38 classes, we expect that a self-dual MCFT $\mathcal{V}_{%
\text{orb}}^{\langle a\rangle }$ with characteristic function $J(\tau )$
exists. Furthermore, we can identify a dual automorphism $a^{*}$ of order $n$
so that $\mathcal{V}^{\Lambda }=((\mathcal{V}^{\Lambda })_{\text{orb}%
}^{\langle a\rangle })_{\text{orb}}^{\langle a^{*}\rangle }$ and where the $%
a^{*}$ centraliser agrees with a corresponding Monster centraliser in all
known cases \cite{T2}. All of this provides evidence that $(\mathcal{V}%
^{\Lambda })_{\text{orb}}^{\langle a\rangle }\simeq \mathcal{V}^{\natural }$
in each construction lending weight to the uniqueness conjecture.

Let us now define the Thompson series $T_{g}(\tau )$ for each $g\in {\bf M}$ 
\begin{eqnarray}
T_{g}(\tau ) &\equiv &\mathrm{Tr}_{\mathcal{H}^{\natural
}}(gq^{L_{0}-1})=Z\left[ 
\begin{array}{c}
g \\ 
1
\end{array}
\right] (\tau )  \label{Tg} \\
&=&\frac{1}{q}+0+[1+\chi _{A}(g)]q+...
\end{eqnarray}
where $\chi _{A}(g)$ is the character of the 196883 dimensional adjoint
representation for ${\bf M}$. The Thompson series for the identity element
is $J(\tau )\,$of (\ref{J}), which is the hauptmodul for the genus zero
modular group \textrm{$SL$}$(2,{\bf Z})$ as already stated.

Conway and Norton \cite{CN} conjectured and Borcherds \cite{B} proved that $%
T_{g}(\tau )$\ is the hauptmodul for some genus zero fixing modular group $%
\Gamma _{g}$. This remarkable property is known as Monstrous Moonshine. In
general, for $g$ of order $n$, $T_{g}(\tau )$ is found to be $\Gamma _{0}(n)$
invariant up to $h^{\text{th}}$ roots of unity where $h$ is an integer with $%
h|n$ and $h|24$ (see Appendix A for the definition of various standard
modular groups). $g$ is a normal element of ${\bf M}$ if and only if $h=1$,
otherwise $g$ is anomalous. $T_{g}(\tau )$ is fixed by some $\Gamma
_{g}\supseteq \Gamma _{0}(N)$ which is contained in the normalizer of $%
\Gamma _{0}(N)$ in $SL(2,{\bf R})$ where $N=nh$ \cite{CN}. This normalizer
contains the Fricke involution $W_{N}$ $:\tau \rightarrow -1/N\tau $. All
classes of ${\bf M}$ can therefore be divided into Fricke and non-Fricke
type according to whether or not $T_{g}(\tau )$ is invariant under the
Fricke involution\textit{. }There are a total of 51 non-Fricke classes of
which 38 are normal and there are a total of 120 Fricke classes of which 82
are normal. We now briefly describe how the genus zero properties of
Monstrous Moonshine can be understood be using the orbifold ideas reviewed
in the last section \cite{T1}, \cite{T2}.

For each of the 38 Leech lattice automorphisms $a$ satisfying the conditions
(i)-(iii) above we can compute the dual automorphism Thompson series $%
T_{a^{*}}$. This agrees precisely with the genus zero series for the 38
non-Fricke normal classes of the Monster where \cite{T2} 
\begin{eqnarray}
T_{a^{*}}(\tau ) &=&\mathrm{Tr}_{\mathcal{H}^{_{\Lambda
}}}(aq^{L_{0}-1})-a_{1} \\
&=&\frac{1}{\eta _{\bar{a}}(\tau )}-a_{1},  \label{Tastar}
\end{eqnarray}
with $\eta _{\bar{a}}(\tau )={\prod_{k|n} }\eta (k\tau )^{a_{k}}$
where $a$ is a lifting of $\overline{a}\in \mathrm{Co}_{0}$ with
characteristic equation $\det (x-\bar{a})=\prod_{k|n}(x^{k}-1)^{a_{k}}$ and
where $n=o(a)=o(\bar{a})$. $\left\{ a_{k}\right\} $ are called the
'Frame-shape' parameters of $\bar{a}$. We can also identify the other 13
non-Fricke classes which are anomalous and find the corresponding correct
genus zero Thompson series \cite{T2}. This is further evidence for the
assertion that $\mathcal{(}\mathcal{V}^{\Lambda }\mathcal{)}_{\text{orb}%
}^{\langle a\rangle }\simeq \mathcal{V}^{\natural }$ implied by the
uniqueness conjecture for $\mathcal{V}^{\natural }$ which we will now assume
from now on.

Consider next $f\in {\bf M}$, a Fricke element of order $n$. For normal
elements we orbifold $\mathcal{V}^{\natural }$ with respect to $\langle
f\rangle $ to obtain a self-dual MCFT $\mathcal{(\mathcal{V}^{\natural })}_{%
\text{orb}}^{\langle f\rangle }$. Assuming $T_{f}(\tau )$ is a hauptmodul
then $\mathcal{(\mathcal{V}^{\natural })}_{\text{orb}}^{\langle f\rangle
}\simeq \mathcal{V}^{\natural }$ for every normal Fricke element \cite{T2}.
The converse is also true, where given that $\mathcal{(\mathcal{V}^{\natural
})}_{\text{orb}}^{\langle f\rangle }\simeq \mathcal{V}^{\natural }$ for some 
$f\in {\bf M}$ then $T_{f}$ is the hauptmodul for a genus zero modular group
containing the Fricke involution \cite{T2}. In general, assuming the
uniqueness conjecture, for all normal elements 
\begin{equation}
\mathcal{V}^{\Lambda } 
\begin{array}{c}
\stackrel{\langle a\rangle }{\rightarrow } \\ 
\stackrel{\langle a^{*}\rangle }{\leftarrow }
\end{array}
\mathcal{V}^{\natural }\stackrel{\langle f\rangle }{\longleftrightarrow }%
\mathcal{V}^{\natural }\Leftrightarrow T_{a^{*}}\text{, }T_{f}\text{ are
hauptmoduls}  \label{Reorb}
\end{equation}
where the arrows represent an orbifolding with respect to the indicated
group.

For any normal Fricke element $f\in {\bf M}$ of order $n$ (\ref{Reorb}) is
equivalent to the following properties for the twisted sector $\mathcal{V}%
_{f}$ :

(i) The $\mathcal{V}_{f}$ vacuum is unique so that $D_{f}=1$ and has
negative vacuum energy $E_{0}^{f}=-1/n$. $\mathcal{V}_{f}$ is then said to
be tachyonic.

(ii) If $f^{r}$ is Fricke then $f^{s}$ is also Fricke where $%
(o(f^{r}),o(f^{s}))=1\,$and $n=o(f^{r})o(f^{s})$.

These conditions are then sufficient to supply all the poles and residues of 
$T_{f}$ so that $T_{f}$ is a hauptmodul for a genus zero fixing group which
includes the Fricke involution \cite{T1,T2}. The genus zero property for an
anomalous class of ${\bf M}$, which corresponds to the Harmonic formula of 
\cite{CN}, is described in \cite{T2}.

The orbifold method can also be employed \cite{T1} to explain the Conway
Norton Power Map\ Formula, which is a property of the Thompson series,
independent of the hauptmodul property:

\textbf{Power Map\ Formula:} Suppose $T_{g}$ is invariant under $\Gamma
_{0}(n|h)+e_{1}$, $e_{2}$,... . Then for any $d$, $T_{g^{d}}$ is invariant
under $\Gamma _{0}(n^{\prime }|h^{\prime })+e_{1}^{\prime }$, $e_{2}^{\prime
}$,... , where $n^{\prime }=n/(n,d)$, $h^{\prime }=h/(h,d)$ and $%
e_{1}^{\prime }$, $e_{2}^{\prime }$,... , are the divisors of $n^{\prime
}/h^{\prime }$ amongst the numbers $e_{1}$, $e_{2}$,... .

\section{Generalised Moonshine}

\subsection{Properties of Generalised Moonshine Functions}

We now consider Generalised Moonshine Functions (GMFs) which are generalised
Thompson series depending on two commuting Monster elements of the form 
\begin{equation}
Z\left[ 
\begin{array}{c}
h \\ 
g
\end{array}
\right] (\tau )=\mathrm{Tr}_{\mathcal{H}_{g}^{\natural }}(hq^{L_{0}-1}),
\label{Z}
\end{equation}
for $h\in C_{g}$ where the action of the lifting of $h$ on $\mathcal{H}%
_{g}^{\natural }$ $\,$is also denoted by $h$ as discussed in section 2.
Norton has conjectured that \cite{N1}:

\noindent{ \bf Generalised Moonshine.} {\it 
The GMF\textbf{\ }(\ref{Z}) is either constant or is a hauptmodul for some
genus zero fixing group\textit{\ }$\Gamma _{h,g}$.
}

\textbf{Note.} We denote $\,$the fixing group for $Z\left[ 
\begin{array}{c}
h \\ 
g
\end{array}
\right] (o(g)\tau )\,$ by $\tilde{\Gamma}_{h,g}$ which is obviously
conjugate to $\Gamma _{h,g}\mathit{.}$

The properties of (\ref{Z}) that follow from the previous section can be
summarised as follows:

(i) When all elements of $\langle g,h\rangle $ $\,$are normal elements of $%
{\bf M}$ then for $\gamma \in SL(2,{\bf Z)}$

\begin{equation}
Z\left[ 
\begin{array}{c}
h \\ 
g
\end{array}
\right] (\gamma \tau )=Z\left[ 
\begin{array}{c}
h^{d}g^{-b} \\ 
h^{-c}g^{a}
\end{array}
\right] (\tau ),  \label{SLT}
\end{equation}
Hence $\Gamma _{h,g}\supseteq \Gamma (o(h),o(g))$ defined in Appendix A. In
practice $\Gamma (o(h),o(g))\triangleleft \Gamma _{h,g}$. In particular,
since $\gamma $ and $-\gamma $ act equally we have

\begin{equation}
Z\left[ 
\begin{array}{c}
h \\ 
g
\end{array}
\right] =Z\left[ 
\begin{array}{c}
h^{-1} \\ 
g^{-1}
\end{array}
\right] ,  \label{Z1}
\end{equation}
which property is known as charge conjugation invariance.

(ii) Given the uniqueness of the twisted sectors for $\mathcal{V}^{\natural
} $, under conjugation by any element $x\in {\bf M}$ then $x(\mathcal{V}%
_{g}^{\natural })x^{-1}$ is isomorphic to $\mathcal{V}_{xgx^{-1}}^{\natural
} $ so that

\begin{eqnarray}
Z\left[ 
\begin{array}{c}
h \\ 
g
\end{array}
\right] &=&\theta (g,h,x)Z\left[ 
\begin{array}{c}
xhx^{-1} \\ 
xgx^{-1}
\end{array}
\right] ,  \label{Z2} \\
\theta (g,h,x) &=&\frac{\phi _{g}(h)}{\phi _{xgx^{-1}}(xhx^{-1})}.
\label{theta}
\end{eqnarray}
with $\phi _{g}$ of (\ref{phi_gh_gahb}) and (\ref{phi_gh_n}).

(iii) For a normal Fricke element $f$ of order $n$ the twisted sector $%
\mathcal{H}_{f}^{\natural }\,$ is unique with vacuum energy $E_{0}^{f}=-1/n$%
. Hence 
\begin{equation}
\phi _{f}(f)=\omega _{n}  \label{FA}
\end{equation}
from (\ref{ghat on twisted vacuum}). Recall that the extension determined by
the phase $\phi _{f}(h)\,$ for each $h\in C_{f}$ is chosen in order to
comply with (\ref{SLT}). If $h$ and $hf$ are in the same $C_{f}$ conjugacy
class where $h=x(hf)x^{-1}$ for some $x\in C_{f}$, then using (\ref{SLT})
for $\gamma =T$, (\ref{Z1}) and (\ref{Z2}) we find that $Z\left[ 
\begin{array}{c}
h \\ 
f
\end{array}
\right] =\theta (f,h,x)Z\left[ 
\begin{array}{c}
hf \\ 
f
\end{array}
\right] $ where $\theta (f,h,x)=\omega _{n}$. In general, conjugating with
respect to elements of $C_{f}$ we find that the GMF is a class function up
to such an $n^{th}$ root of unity i.e. $\theta (g,h,x)\in \langle \omega
_{n}\rangle $ for all $x\in C_{f}$. If, on the other hand, $h$ and $hf$ are
not in the same $C_{f}$ conjugacy class then $\phi _{f}\,$ can be chosen so
that 
\begin{equation}
\phi _{f}(xhx^{-1})=\phi _{f}(h)  \label{phi_class_fun}
\end{equation}
for all $x\in C_{f}$.

From (\ref{Z(h,g) char expansion}) for $h\in C_{f}$ we have 
\begin{eqnarray}
Z\left[ 
\begin{array}{c}
h \\ 
f
\end{array}
\right] (\tau ) &=&q^{-1/n}\sum_{m=0}^\infty\chi
_{f,m}^{(1)}(h)q^{m}+\sum_{j=0,j\neq 1}^{n-1}q^{-j/n}%
\sum_{m=1}^\infty \chi _{f,m}^{(j)}(h)q^{m} 
\nonumber \\
&=&\phi _{f}(h)[q^{-1/n}+0+\sum{s=1}^\infty
a_{f,s}(h)q^{s/n}].  \label{Z(h,Fricke)}
\end{eqnarray}
The coefficient $a_{f,s}(h)=\chi _{f,m}^{(j)}(h)/\phi _{f}(h)$ is called a
head character for the given GMF where $s=mn-j$. Then $a_{f,s}(h)$ is a
character for $G_{f}$ where $G_{f}=C_{f}/\langle f\rangle \,$ for $s=-1\
{\rm mod}\ n$ and otherwise is a character for $C_{f}$. For example, if we
choose $f=2+$ (the $2A\,$ Monster element) then $C_{f}=2.B$ and $G_{f}=B$,
the Baby Monster. Then $a_{f,s}(h)$ is a character for $B$ for odd $s$ and
is a character for $2.B$ for even $s$ \cite{DLM2}. This implies that the
decomposition of (\ref{Z(h,Fricke)}) into the irreducible characters of $%
C_{g}$ involves only those characters obeying such conditions. This
observation is confirmed in Appendix B for $f=2+$ and $3+$ where explicit
character expansions are given.

(iv) For a normal non-Fricke element $a^{*}$ of order $n$ we find the $S$
transformation of $T_{a^{*}}$ of (\ref{Tastar}) results in 
\begin{equation}
Z\left[ 
\begin{array}{c}
1 \\ 
a^{*}
\end{array}
\right] (\tau )=-a_{1}+O(q^{1/n}).  \label{Seta}
\end{equation}
i.e. $\mathcal{H}_{a^{*}}^{\natural }$ has vacuum energy $E_{0}^{a^{*}}=0.$
Furthermore, from (\ref{ghat on twisted vacuum}) 
\begin{equation}
\rho _{a^{*}}^{0}(a^{*})=\mathbf{1}  \label{Unit}
\end{equation}
with vacuum degeneracy $D_{a^{*}}=\dim (\rho _{a^{*}}^{0})=-a_{1}$.
Furthermore if $(a^{*})^{B}$ is Fricke for some $B$ then from (\ref{Tastar})
we find that $o(a^{*})||n$ i.e. $o(a^{*})|n$ and $(n,o(a^{*}))=1$.

(v) The value of (\ref{Z}) at any parabolic cusp $a/c$ with $(a,c)=1$ is
determined by the vacuum energy of the $g^{a}h^{-c}$ twisted sector from
(i). If $g^{a}h^{-c}\,$is Fricke then (\ref{Z}) is singular at $a/c$ \cite
{N1} with residue $\phi _{g^{a}h^{-c}}(h^{d}g^{-b})$. (\ref{Z}) is
holomorphic at all other points on ${\bf H}$. Once these singularities are
known, then (\ref{Z}) can be analysed to check whether it is constant or is
a hauptmodul for an appropriate genus zero modular group. If all elements of 
$\langle g,h\rangle $ are Non-Fricke then there are no singular cusps so
that $Z\left[ 
\begin{array}{c}
h \\ 
g
\end{array}
\right] $ is holomorphic on ${\bf H}/\Gamma (o(h),o(g))$ and hence is
constant. This accounts for the constant GMFs referred to in the Generalised
Moonshine Conjecture above. We therefore assume from now on that at least
one element of $\langle g,h\rangle $ $\,$is Fricke which we chose to be $g$
without loss of generality. The singularities and residues of $Z\left[ 
\begin{array}{c}
h \\ 
g
\end{array}
\right] $ are then constrained by certain orbifolding constraints which we
discuss below in terms of two consistency theorems.

(vi) The operators $\{L_{0},L_{-1},L_{1}\}\,$generate an $sl(2,{\bf C})$
\thinspace sub-algebra acting on $\mathcal{H}_{g}=\bigoplus_{j=0}^{n-1}%
\bigoplus_{m=0}^{\infty }\mathcal{H}_{g,m}^{(j)}$ where $L_{-1}(\mathcal{H}%
_{g,m}^{(j)})\subseteq \mathcal{H}_{g,m+1}^{(j)}$ and $L_{1}(\mathcal{H}%
_{g,m}^{(j)})\subseteq \mathcal{H}_{g,m-1}^{(j)}$ so that 
\begin{equation}
\mathcal{H}_{g,m}^{(j)}=\ker _{\mathcal{H}_{g,m}^{(j)}}L_{1}\oplus \mathrm{im%
}_{\mathcal{H}_{g,m-1}^{(j)}}L_{-1}.
\end{equation}
Hence the corresponding representations for $\mathrm{Aut}(\mathcal{V}_{g})$
are related where $\rho _{g,m}^{(j)}=\rho _{g,m-1}^{(j)}+\tilde{\rho}%
_{g,m}^{(j)}$ where $\tilde{\rho}_{g,m}^{(j)}$ is the representation formed
by $\ker _{\mathcal{H}_{g,m}^{(j)}}L_{1}$. Therefore every character $\chi
_{g,m}^{(j)}(h)$ for $h\in \mathrm{Aut}(\mathcal{V}_{g})$ obeys the property 
\begin{equation}
\chi _{g,m}^{(j)}(h)=\chi _{g,m-1}^{(j)}(h)+\tilde{\chi}_{g,m}^{(j)}(h).
\label{chim}
\end{equation}
In Appendix B we observe this property for the first ten head characters $%
a_{g,mn-j}(h)=\chi _{g,m}^{(j)}(h)/\phi _{g}(h)$ for $g=p+$ and $p=2$ and $3$%
. For $p=5$ and $7$ this property can be observed in \cite{Q1},\cite{Q2}.

\subsection{Two Consistency Theorems}

In this section we prove two theorems which simplify our analysis of GMFs.
Their purpose is to identify the phases $\phi _{h^{-c}g^{a}}(h^{d}g^{-b})\,$
for $h^{-c}g^{a}\,$ Fricke appearing in (\ref{SLT}). This will be achieved
by considering the consistency of orbifolding $\mathcal{V}^{\natural }\,$%
with respect to $\langle g,h\rangle \,$ under various choices of independent
generators.

\begin{thm}
\label{theorem_proj1} Let $g,h\in {\bf M}$ be independent commuting elements
where both $g$ and $h$ are Fricke such that $\phi _{g}(h)=1$ and where all
elements of $\langle g,h\rangle \,$ are normal. Let $u$, $v$ be any
independent generators for $\langle g,h\rangle $. Then:
\end{thm}

\begin{description}
\item[(i)]  $\mathcal{V}_{\text{orb}}^{\langle g,h\rangle }\simeq \mathcal{V}%
^{\natural }.$

\item[(ii)]  If $u$ is non-Fricke then $\phi _{u^{A}v^{B}}(u)\neq 1$ for all
Fricke elements $u^{A}v^{B}$ except possibly when $o(v^{B})||o(v)$ with $%
(B,o(v))\neq 1$.

\item[(iii)]  If $u$ is Fricke then there is a unique $A$ mod $o(u)$ such
that $u^{A}v$ is Fricke with $\phi _{u^{A}v}(u)=1$ and $o(u^{A}v)=o(v)$.
\end{description}

\textbf{Proof }(i) Since $h$ is Fricke and normal we have $\mathcal{V}_{%
\text{orb}}^{\langle h\rangle }\simeq \mathcal{V}^{\natural }$ from (\ref
{Reorb}). By assumption $\phi _{g}(h)=1\,$and hence $(\mathcal{V}_{\text{orb}%
}^{\langle h\rangle })_{g}=\mathcal{P}_{h}(\mathcal{V}_{g}^{\natural }\oplus 
\mathcal{V}_{gh}^{\natural }\oplus \mathcal{V}_{gh^{2}}^{\natural }\oplus
...)$ is tachyonic i.e. has negative vacuum energy. Therefore $g$ acts as a
normal Fricke element on $\mathcal{V}_{\text{orb}}^{\langle h\rangle }\simeq 
\mathcal{V}^{\natural }$ and orbifolding the latter with respect to $g$ we
find that $\mathcal{V}_{\text{orb}}^{\langle g,h\rangle }\simeq \mathcal{V}%
^{\natural }$ from (\ref{Reorb}) again.

(ii) Since $u,v\,$ are independent $\mathcal{V}_{\text{orb}}^{\langle
u,v\rangle }=(\mathcal{V}_{\text{orb}}^{\langle u\rangle })_{\text{orb }%
}^{\langle v\rangle }\simeq \mathcal{V}^{\natural }$ $\,$from (i). With $u$
Non-Fricke then $\mathcal{V}_{\text{orb}}^{\langle u\rangle }\simeq \mathcal{%
V}^{\Lambda }$ and hence $(\mathcal{V}_{\text{orb}}^{\langle u\rangle
})_{v^{B}}=\mathcal{P}_{u}(\mathcal{V}_{v^{B}}^{\natural }\oplus \mathcal{V}%
_{uv^{B}}^{\natural }\oplus \mathcal{V}_{u^{2}v^{B}}^{\natural }\oplus ...)$
is non-tachyonic except possibly when $o(v^{B})||o(v)$ and $(B,o(v))\neq 1$.
But if $u^{A}v^{B}$ is Fricke, then $\mathcal{V}_{u^{A}v^{B}}^{\natural }$
is tachyonic so that $\mathcal{P}_{u}=0$ on the corresponding vacuum sector
i.e. $\phi _{u^{A}v^{B}}(u)\neq 1$ except possibly when $o(v^{B})||o(v)$ and 
$(B,o(v))\neq 1$.

(iii) With $u$ Fricke then $\mathcal{V}_{\text{orb}}^{\langle u\rangle
}\simeq \mathcal{V}^{\natural }.$ But $(\mathcal{V}_{\text{orb}}^{\langle
u\rangle })_{\text{orb }}^{\langle v\rangle }\simeq \mathcal{V}^{\natural }$
implies that $v$ acts a Fricke element of order $o(v)$ on $\mathcal{V}_{%
\text{orb}}^{\langle u\rangle }$ since $u,v\,$ are independent. Hence $(%
\mathcal{V}_{\text{orb}}^{\langle u\rangle })_{v}=\mathcal{P}_{u}(\mathcal{V}%
_{v}^{\natural }\oplus \mathcal{V}_{uv}^{\natural }\oplus \mathcal{V}%
_{u^{2}v}^{\natural }\oplus ...)$ is tachyonic which is possible iff $%
\mathcal{P}_{u}=1$ on precisely one of the tachyonic vacuum sectors $%
\mathcal{V}_{u^{A}v}^{(0)}$ for some unique $A$ mod $o(u)$ where $%
o(v)=o(u^{A}v).$ $\qed $

Example. Consider the orbifolding of $\mathcal{V}^{\natural }\,$with respect
to $\langle f\rangle $ for $f$ Fricke of non-prime order $n$. Let $r||n$ and
consider $u=f^{r}$ of order $s=n/r$ and $v=f^{s}$ of order $r$. Then $u,v\ $%
are independent generators of $\langle f\rangle $ and Theorem \ref
{theorem_proj1} (iii) with $g=f$ and $h=1$ implies that if $u$ is Fricke
then there is a unique $A\ {\rm mod}\ s$ such that $u^{A}v=f^{Ar+s}$ is
Fricke of order $r$. But $(r,s)=1$ implies that $A=0\ {\rm mod}\ s$ so that 
$v$ is Fricke. This is the Atkin-Lehner closure property for Thompson series 
\cite{T1} i.e. if $f$ and $f^{r}$ are Fricke for $r||n\,$then $f^{s}$ is
also Fricke for $s=n/r$.

\begin{thm}
\label{theorem_proj2} Let $g,h\in {\bf M}$ be independent commuting elements
where $g$ is Fricke and $h$ is Non-Fricke such that $\phi _{g}(h)=1$ and
where all elements of $\langle g,h\rangle \,$ are normal. Let $u$, $v$ be
any independent generators for $\langle g,h\rangle $. Then:
\end{thm}

\begin{description}
\item[(i)]  $\mathcal{V}_{\text{orb}}^{\langle g,h\rangle }\simeq \mathcal{V}%
^{\Lambda }$.

\item[(ii)]  If $u$ is Fricke then $\phi _{u^{A}v^{B}}(u)\neq 1$ for all $%
u^{A}v^{B}$ Fricke except possibly when $o(v^{B})||o(v)$ with $(B,o(v))\neq 1
$.
\end{description}

\textbf{Proof }(i) Since $g$ is Fricke $\mathcal{V}_{\text{orb}}^{\langle
g\rangle }\simeq \mathcal{V}^{\natural }$ from (\ref{Reorb}). Then either $%
\mathcal{V}_{\text{orb}}^{\langle g,h\rangle }=\mathcal{P}_{h}(\mathcal{V}_{%
\text{orb}}^{\langle g\rangle }\oplus (\mathcal{V}_{\text{orb}}^{\langle
g\rangle })_{h}\oplus ...)$ $\simeq \mathcal{V}^{\natural }$ or $\mathcal{V}%
_{\text{orb}}^{\langle g,h\rangle }\simeq \mathcal{V}^{\Lambda }$ when $h$
acts as Fricke or Non-Fricke element on $\mathcal{V}_{\text{orb}}^{\langle
g\rangle }$ from (\ref{Reorb}) again. Assume that $\mathcal{V}_{\text{orb}%
}^{\langle g,h\rangle }\simeq \mathcal{V}^{\natural }\,$ and consider the
alternative composition of orbifoldings $\mathcal{V}_{\text{orb}}^{\langle
g,h\rangle }=\mathcal{P}_{g}(\mathcal{V}_{\text{orb}}^{\langle h\rangle
}\oplus (\mathcal{V}_{\text{orb}}^{\langle h\rangle })_{g}\oplus ...)$ where
since $h$ is non-Fricke, $\mathcal{V}_{\text{orb}}^{\langle h\rangle }\simeq 
\mathcal{V}^{\Lambda }$. But the condition $\phi _{g}(h)=1$ implies that the
vacuum energy of $(\mathcal{V}_{\text{orb}}^{\langle h\rangle })_{g}=%
\mathcal{P}_{h}\mathcal{V}_{g}\oplus ...$ is negative which is impossible
according to (\ref{E0a}). Therefore $\mathcal{V}_{\text{orb}}^{\langle
g,h\rangle }\simeq \mathcal{V}^{\Lambda }$.

(ii) For independent $u,v$, $\mathcal{V}_{\text{orb}}^{\langle u,v\rangle }=(%
\mathcal{V}_{\text{orb}}^{\langle u\rangle })_{\text{orb }}^{\langle
v\rangle }\simeq \mathcal{V}^{\Lambda }$ from (i) where $u$ Fricke implies $%
\mathcal{V}_{\text{orb}}^{\langle u\rangle }\simeq \mathcal{V}^{\natural }$.
Hence $v$ is Non-Fricke in its action on $\mathcal{V}_{\text{orb}}^{\langle
u\rangle }$ and so $(\mathcal{V}_{\text{orb}}^{\langle u\rangle })_{v^{B}}=%
\mathcal{P}_{u}(\mathcal{V}_{v^{B}}^{\natural }\oplus \mathcal{V}%
_{uv^{B}}^{\natural }\oplus \mathcal{V}_{u^{2}v^{B}}^{\natural }\oplus ...)$
is non-tachyonic except possibly when $o(v^{B})||o(v)$ and $(B,o(v))\neq 1$.
But if $u^{A}v^{B}$ is Fricke, then $\mathcal{V}_{u^{A}v^{B}}^{\natural }$
is tachyonic so that $\mathcal{P}_{u}=0$ on the corresponding vacuum sector
i.e. $\phi _{u^{A}v^{B}}(u)\neq 1$ except possibly when $o(v^{B})||o(v)$ and 
$(B,o(v))\neq 1$. $\qed $

Note. There is no corresponding statement to (iii) in Theorem \ref
{theorem_proj1}. This is because not all self-orbifoldings of the Leech
theory are of 'Fricke' type e.g. the involution of $\mathrm{Co}_{0}$ with
frame shape $2^{16}/1^{8}\,$has a lifting $a\in \mathrm{Aut}(\mathcal{V}%
^{\Lambda })\,$ with non-tachyonic vacuum energy $E_{a}^{0}=0$ and vacuum
degeneracy $D_{a}=16\,$ whereas $(\mathcal{V}^{\Lambda })_{\text{orb}%
}^{\langle a\rangle }\simeq \mathcal{V}^{\Lambda }$.

\subsection{Symmetries of Rational GMFs for $g=p+$}

We now consider GMFs (\ref{Z}) where $g$ is assumed to be Fricke and of
prime order $p$ i.e. $g=p+$ in Conway-Norton notation (see Appendix A) for $%
p=2,3,5,...,31,41,47,59,71$. If $\langle g,h\rangle =\langle u\rangle $ for
some $u\in {\bf M}$ then (\ref{Z}) can always be transformed to a regular
Thompson series (\ref{Tg}) via a modular transformation (\ref{SLT}) as
follows

\begin{equation}
Z\left[ 
\begin{array}{c}
h \\ 
g
\end{array}
\right] (\tau )=Z\left[ 
\begin{array}{c}
u^{A} \\ 
1
\end{array}
\right] (\gamma \tau ),\quad  \label{<g,h>=<u>}
\end{equation}
where $g=u^{-C}$, $h=u^{D}$ for some $C,D$ $\,$and where $A=(C,D)$, $c=C/A$, 
$d=D/A\,$ and $a,b$ are chosen so that $ad-bc=1$. In particular, this is
always possible when $p\,$and $o(h)$ are coprime. An example of this
phenomenon is discussed later on in Section 4.4. In these cases, the genus
zero property for the GMF therefore follows from that for a regular Thompson
series (\ref{Tg}). The GMFs with non-trivial genus zero behaviour then occur
for $h\in C_{g}$ where $o(h)=pk$ for $k\geq 1$. These only exist for $p\leq
13$. Furthermore, from now on we restrict our analysis to $k=1\,$ or $k$
prime alone.

We now make the further restriction of considering rational GMFs i.e. those
with $q$ expansion (\ref{Z(h,Fricke)}) with rational (and therefore
integral) coefficients. No rational GMFs for $o(h)=pk$ occur for $p=11$ and $%
13$ from explicit calculations of the head characters. We therefore only
consider $p=2,3,5$ and $7\,$ from now on. This firstly implies that the GMF
for $p>2$ has leading $q$ expansion in normalised function form \cite{FMN} 
\begin{equation}
Z\left[ 
\begin{array}{c}
h \\ 
g
\end{array}
\right] (\tau )=\frac{1}{q^{1/p}}+0+O(q^{1/p}),  \label{Z(h,g) =1/q^(1/p)}
\end{equation}
i.e. $\phi _{g}(h)=1$. We may also assume this normalised function form for $%
p=2$ by relabelling as discussed below. This choice is also sufficient to
ensure that $g,h$ are independent generators of $\langle g,h\rangle $ since
otherwise $h^{k}=g^{A}$ for some $A\neq 0$ ${\rm mod}$ $p$ which implies
the contradictory relation $(\phi _{g}(h))^{k}=(\phi _{g}(g))^{A}=\omega
_{p}^{A}$ from (\ref{phi_gh_gahb}). Then the centraliser $C_{g}$ and
sporadic finite simple group $G_{g}=C_{g}/\langle g\rangle \,$ are as
follows:

\begin{tabular}{|l|l|l|l|}
\hline
$g=p+$ & $C_{g}$ & $G_{g}$ & Name \\ \hline
$2+$ & $\mathrm{2.B}$ & $\mathrm{B}$ & Baby Monster \\ \hline
$3+$ & $\mathrm{3.Fi}$ & $\mathrm{Fi}$ & Fischer \\ \hline
$5+$ & $\mathrm{5\times HN}$ & $\mathrm{HN}$ & Harada-Norton \\ \hline
$7+$ & $\mathrm{7\times He}$ & $\mathrm{He}$ & Held \\ \hline
\end{tabular}

Table 1. $g=p+$ centralisers for $p=2,3,5,7$.

Each element $l\in G_{g}$ is the image of $p$ distinct elements, which we
denote by $(g^{a},l)\in C_{g}$ for $a=1,...p$ where $(g^{a_{1}},l_{1}).$ $%
(g^{a_{2}},l_{2})=(g^{a_{2}+a_{2}+c(l_{1},l_{2})},l_{1}l_{2})$ with
non-trivial cocycle $c(l_{1},l_{2})$ for $p=2,3$. If $(g^{a_{1}},l_{1})%
\stackrel{C_{g}}{\sim }$ $(g^{a_{2}},l_{2})$ then clearly $l_{1}\stackrel{%
G_{g}}{\sim }l_{2}$ and hence the conjugacy classes of $C_{g}$ can be
labelled (non-uniquely) by the conjugacy classes of $G_{g}$ as in the ATLAS 
\cite{CCNPW}. We define the $(1,G_{g})$ classes of $C_{g}$ to be those
classes displayed in the ATLAS. For $p=2,3$ where $C_{g}$ is a non-trivial
extension of $G_{g}$ then $h\stackrel{C_{g}}{\sim }$ $gh$ for some elements $%
h\in C_{g}$ and the corresponding $(1,G_{g})$ class contains elements of the
form $(g^{a},l)$ for all $a$. Any $G_{g}$ irreducible character $\chi
^{G_{g}}$ becomes a $C_{g}$ irreducible character with $\chi (g^{a},l)=\chi
^{G_{g}}(l)\,$ for $a=1,...p$ whereas for the remaining $C_{g}$ irreducible
characters we have $\chi ^{C_{g}}(g^{a},l)=\omega _{p}^{a}\chi ^{C_{g}}(1,l)$%
. The irreducible characters for the $(1,G_{g})$ classes are given
unambiguously in the ATLAS for $p>2$ \cite{CCNPW}. For $p=2$, we define the $%
(1,G_{g})\,$ classes to be those classes with characters as given in the
ATLAS for $2.\mathrm{B}$.

We will now restrict our attention to those elements $h$ $\,$where $o(h)=pk\,
$ for $k=1\,$ and $k$ prime and $h$ an element of a $(1,G_{g})$ class of $%
C_{g}$. We will find all the corresponding GMFs below using orbifold
considerations and demonstrate that they must all have genus zero fixing
groups. Towards this aim we firstly prove a theorem concerning the phase
multipliers $\phi _{g}(h)$.

\begin{thm}
\label{theorem_phi} For $g=p+$, $p=2,3,5\,$and $7$ and $o(h)=pk$, $k=1\,$ or 
$k$ prime and assuming (\ref{phi_gh_gahb}) then for $h$ normal $\phi
_{g}(h)\in \langle \omega _{p}\rangle $.
\end{thm}

\textbf{Proof.} If $k=1$ then the result is obvious from (\ref{phi_gh_gahb})
since $\phi _{g}(h)^{p}=\phi _{g}(h^{p})=1$. It is necessary to separately
consider (a) $k=p\,$and (b) $k\neq p$, $k>1$.

(a) Assume $k=p$. For $p=2$ using charge conjugation (\ref{Z1}) we have $%
\phi _{g}(h)=\phi _{g}(h^{-1})\,$ since $g^{-1}=g$ so that $\phi _{g}(h)\in
\langle \pm 1\rangle $ for all normal $h$. However, note that with $h=4|2+$,
an anomalous Monster class, we have $g=h^{2}=2+$ and $h$ projects down to
the $2C$ class of $\mathrm{B}$. However in this case $\phi _{g}(h)^{2}=\phi
_{g}(h^{2})=\phi _{g}(g)=-1$ and hence $\phi _{g}(h)=\pm i$ so that (\ref{Z1}%
) doesn't hold. This anomalous case is discussed later on in Section 4.5.

For $p>2$, $H=h^{p}$ is of order $p$ and hence $\phi _{g}(H)\in \langle
\omega _{p}\rangle $. From (\ref{phi_gh_gahb}) we have $\phi _{g}(H)=\phi
_{g}(h)^{p}$. On the other hand $\phi _{g}(xHx^{-1})=\phi _{g}(xhx^{-1})^{p}$
for all $x\in C_{g}$. But then $\theta (g,H,x)=\theta (g,h,x)^{p}=1$ and so
from (\ref{theta}) $\phi _{g}(H)=$ $\phi _{g}(xHx^{-1})$. For $p>2$ we also
find that $H$ must be an element of a class of type $(1,G_{g})$. For $p=5\,$
and $7$ this is obvious since the extension of $G_{g}$ to $C_{g}$ is
trivial. For $p=3$, all classes of order $3\,\,$in $\mathrm{Fi\,}$ lift to
classes of order $3$ in $3.\mathrm{Fi}$ from the ATLAS \cite{CCNPW} and
hence the same result follows. From the ATLAS it is also clear that every
class of order $p$ and of type $(1,G_{g})$ is conjugate to at least one its
powers. Hence $\phi _{g}(H)=\phi _{g}(H^{a})\,$ for some $a\neq 1\ {\rm mod}\ p$%
. But since $\phi _{g}(H)\in \langle \omega _{p}\rangle $ it follows that $%
\phi _{g}(H)=\phi _{g}(H)^{a}=1$ for $p>2$. Hence $\phi _{g}(h)\in \langle
\omega _{p}\rangle $ as claimed.

(b) Lastly consider $k\neq p$, $k>1$. Then $H=h^{p}$ is of order $k$ coprime
to $p$ with $\phi _{g}(H)\in \langle \omega _{k}\rangle $ since $\phi
_{g}(H)^{k}=\phi _{g}(H^{k})=1.$ But in this case $H$ and $g$ are coprime
and hence the GMF $Z\left[ 
\begin{array}{c}
H \\ 
g
\end{array}
\right] \,$ can be directly found from the Thompson series for a normal
class of order $pk$ i.e. either (i) $H=k+\,$with $gH=pk+$ or (ii) $H=k-$
with $gH=pk+p\,$ from the power map formula.

(i) $H=k+$ and $gH=pk+$. Since $gH$ is Fricke we have $\mathcal{V}_{\text{orb%
}}^{\langle g,H\rangle }=\mathcal{V}_{\text{orb}}^{\langle gH\rangle }\simeq 
\mathcal{V}^{\natural }$. Consider $(\mathcal{V}_{\text{orb}}^{\langle
H\rangle })_{g}=\mathcal{P}_{H}(\mathcal{V}_{g}^{\natural }\oplus \mathcal{V}%
_{gH}^{\natural }\oplus ...\oplus \mathcal{V}_{gH^{k-1}}^{\natural })$ which
must be a tachyonic Fricke twisted sector of $\mathcal{V}_{\text{orb}%
}^{\langle H\rangle }\,$of order $p$ since $\mathcal{V}_{\text{orb}%
}^{\langle H\rangle }\simeq \mathcal{V}^{\natural }$. However, $gH^{a}$ is
of order $p$ for $a=0\ {\rm mod}\ p$ only so that $\mathcal{P}_{H}(\mathcal{V}%
_{g}^{\natural })$ $\,$is tachyonic and hence $\phi _{g}(H)=1$.

(ii) $H=k-$ and $gH=pk+p$. Since $gH$ is Non-Fricke we have $\mathcal{V}_{%
\text{orb}}^{\langle gH\rangle }\simeq \mathcal{V}^{\Lambda }$. Hence $%
Z\left[ 
\begin{array}{c}
gH \\ 
1
\end{array}
\right] (\tau )$ $\,$is the known Thompson series 
\begin{equation}
T_{pk+p}(\tau )=\ \left[ \frac{\eta (\tau )\eta (p\tau )}{\eta (k\tau )\eta
(pk\tau )}\right] ^{\frac{24}{(k-1)(p+1)}}+\frac{24}{(k-1)(p+1)},
\end{equation}
from (\ref{Tastar}). This is invariant under the Atkin-Lehner transformation 
$W_{p}$ (A.5) of Appendix A which we can choose as $W_{p}=\left( 
\begin{array}{ll}
p & 1 \\ 
ckp & dp
\end{array}
\right) $ where $dp-ck=1$ for coprime $p,k$. Then 
\begin{eqnarray}
Z\left[ 
\begin{array}{c}
gH \\ 
1
\end{array}
\right] (\tau ) &=&Z\left[ 
\begin{array}{c}
gH \\ 
1
\end{array}
\right] (W_{p}(\tau ))  \nonumber \\
&=&Z\left[ 
\begin{array}{c}
H \\ 
g
\end{array}
\right] (p\tau )=\frac{1}{q}+0+...
\end{eqnarray}
so that $\phi _{g}(H)=1$ and hence $\phi _{g}(h)\in \langle \omega
_{p}\rangle $. Note that this argument can also be directly used in (i) if
we consider the $W_{p}$ invariance for $T_{pk+}$.$\qed $

\begin{cor}
\label{corollary_phi1} With $g,h$ as above in Theorem \ref{theorem_phi} then

\begin{enumerate}
\item[(a)]  If $h$ and $gh$ are conjugate in $C_{g}$ we can choose $\phi
_{g}(h)=1$ for $p=2,3$.

\item[(b)]  If $h$ and $gh$ are not conjugate in $C_{g}$ then for $p=3,5,7$
we have $\phi _{g}(h)=1$ iff $h$ is an element of a $(1,G_{g})$ class of $%
C_{g}$.
\end{enumerate}
\end{cor}

\textbf{Proof.} (a) As noted in Section 3.1 (iii) above if $h$ and $gh$ are
conjugate in $C_{g}$ (and hence both belong to a $(1,G_{g})$ class) then $%
Z\left[ 
\begin{array}{c}
gh \\ 
g
\end{array}
\right] =\omega _{p}Z\left[ 
\begin{array}{c}
h \\ 
g
\end{array}
\right] $ and similarly for any $g^{a}h$ for all $a$. Hence from Theorem \ref
{theorem_phi} $\phi _{g}(g^{a}h)=1\,$ for some $a$ so that after the
relabeling $g^{a}h\rightarrow h$ we choose $\phi _{g}(h)=1$.

(b) If $h$ and $gh$ are not conjugate in $C_{g}$ then $h$ is conjugate to $%
h^{a}$ for some $a\neq 1\ {\rm mod}\ pk$ iff $h$ is a class of type $%
(1,G_{g})$ from the Atlas \cite{CCNPW}. Furthermore from (\ref{phi_class_fun}%
) for $p>2$ it follows that $\phi _{g}(h)=\phi _{g}(h^{a})\,$ for some $%
a\neq 1\ {\rm mod}\ pk$ iff $\phi _{g}(h)=1$ and hence the result follows.$%
\qed $

\textbf{Remark}. For $p=2$ if $h$ and $gh$ are not conjugate in $C_{g}$ then
we may choose $\phi _{g}(h)=1$ in conjunction with the definition of the $%
(1,G_{g})\,$ classes for $2.\mathrm{B}$.

\begin{cor}
\label{corollary_phi2} For $p=3,5,7$ then $Z\left[ 
\begin{array}{c}
h \\ 
g
\end{array}
\right] \,$rational implies $\phi _{g}(h)=1\,$and $h$ is a member of a class
of type $(1,G_{g})$.
\end{cor}

\textbf{Proof. }From Theorem \ref{theorem_phi} $\phi _{g}(h)\in \langle
\omega _{p}\rangle $ which is rational for $\phi _{g}(h)=1$ only for $p>2$.
Hence from Corollary \ref{corollary_phi1} we have $h$ is a member of a class
of type $(1,G_{g})$.$\qed $

We make the following useful observations concerning the irreducible
characters for the groups of Table 1 which can be checked by inspecting the
appropriate ATLAS character tables \cite{CCNPW}.

\begin{lem}
\label{lemma_Cg_characters} The irreducible characters $\chi $ for the
groups $C_{p+}$ of Table 1 enjoy the following properties for $(1,G_{g})$
classes:

\begin{description}
\item[(a)]  Any irrationality for $\chi $ is quadratic i.e. for $h\in C_{p+}$%
, $\chi (h)=a\pm\sqrt{b}$ for some $a,b\in {\bf Q}$ where $b=0$ for rational $%
\chi (h)$.

\item[(b)]  The number, $N_{\chi }$, of independent irreducible characters
of a given dimension is $N_{\chi }=1$ or $2$ or in the case of the Held
group possibly $N_{\chi }=3$. If $N_{\chi }=2$ then the two characters, $%
\chi $ and $\bar{\chi}$ (say), are irrational and algebraically conjugate
i.e. $\bar{\chi}(h)\equiv a\mp \sqrt{b}$ for all $h$ in the notation of (a). In
the case of the Held group there are three characters of dimension $1275$
where two of the characters are irrational and algebraically conjugate and
the other is rational.
\end{description}

\textbf{Note.} A character is said to be irrational if it is irrational for
at least one conjugacy class. Otherwise, it is said to be rational. Clearly
if $N_{\chi }=1$ then the irreducible character $\chi $ is rational since $%
\bar{\chi}\,$ is also an irreducible character.
\end{lem}

\textbf{Examples.} The Fischer group has two irreducible algebraically
conjugate characters $\chi _{6},\chi _{7}$ of dimension $1603525$ and two
inequivalent conjugacy classes $23A$, $23B$ \cite{CCNPW} such that $\chi
_{6}(23A)=\chi _{7}(23B)=(-1+i\sqrt{23})/2$ and $\chi _{6}(23B)=\chi
_{7}(23A)=(-1-i\sqrt{23})/2$.

\begin{thm}
\label{theorem_conjugation} Consider $g=p+\ $for $p=3,5,7$ and $o(h)=pk$, $%
k=1$ or $k$ prime. Suppose that $xC_{g}x^{-1}=C_{g}\,$ for some $x\in {\bf M}
$ and that the GMF $Z\left[ 
\begin{array}{c}
h \\ 
g
\end{array}
\right] (\tau )$ is rational. Then 
\begin{equation}
Z\left[ 
\begin{array}{c}
h \\ 
g
\end{array}
\right] =Z\left[ 
\begin{array}{c}
xhx^{-1} \\ 
g
\end{array}
\right] .  \label{general conjugation}
\end{equation}
\end{thm}

\smallskip \textbf{Proof.} Since the GMF is rational from Corollary \ref
{corollary_phi2} then $\phi _{g}(h)=1$ \thinspace and $h$ is a member of a $%
(1,G_{g})\,$ class. Furthermore $h$ is conjugate to $h^{a}$ for some $a\neq
1\ {\rm mod}\ pk$ from an inspection of the ATLAS character tables \cite
{CCNPW}. Hence $xhx^{-1}$ is conjugate to $(xhx^{-1})^{a}$ and is also
member of a $(1,G_{g})\,$ class of $C_{g}$. Therefore $\phi
_{g}(xhx^{-1})=1\,$also. For any irreducible character $\chi $ of $C_{g}$
and for all $h\in C_{g}$, $\chi ^{x}(h)\equiv \chi (xhx^{-1})$ $\,$is
another irreducible character of the same dimension. Furthermore $%
xC_{g}x^{-1}=C_{g}$ implies that the number of elements in the conjugacy
classes for $h$ and $xhx^{-1}$ in $C_{g}$ are equal. We now show that either 
$\chi ^{x}=$ $\chi $ or $\bar{\chi}$, $\,$the algebraic conjugate. Suppose
that $\chi ^{x}\neq \chi $. From Lemma \ref{lemma_Cg_characters} then either
(i) they are both irrational and algebraically conjugate with $\chi ^{x}=$ $%
\bar{\chi}$ or else (ii) $G_{g}=\mathrm{He}\,$ and $\chi \,$ is of dimension
1275 where one is irrational $\chi $, say, and the other $\chi ^{x}$ is the
unique rational character of dimension 1275. However the character table 
\cite{CCNPW} for $\mathrm{He}$ reveals that (ii) cannot be true since there
are no two conjugacy classes of $\mathrm{He}$ with the same number of
elements for which a 1275 dimensional character $\chi $ is irrational and
the other $\chi ^{x}$ is rational. Hence $\chi ^{x}$ $=\bar{\chi}\,$ if $%
\chi ^{x}\neq \chi $.

The head character $a_{g,mp-j}(h)$ coefficient of $q^{m-j/p}\,$ in (\ref
{Z(h,Fricke)}) is rational by assumption and is some integer linear
combination of the irreducible characters for $C_{g}.$ If a given
irreducible character $\chi (h)$ is irrational then $\chi \,$ and $\bar{\chi}%
\,$ must appear with the same (possibly zero) multiplicity and $\chi +\bar{%
\chi}=\chi ^{x}+(\bar{\chi})^{x}$ whereas if $\chi (h)$ is rational then $%
\chi (h)=\chi ^{x}(h)\,$. It is therefore clear that $%
a_{g,mp-j}(h)=a_{g,mp-j}(xhx^{-1})$ for each head character and hence the
result follows. $\qed $

\begin{thm}
\label{theorem_GMF_power_conjugation} Consider $g=p+\ $for $p=2,3,5,7$ and $%
o(h)=pk$, $k=1$ or $k$ prime. If the GMF $Z\left[ 
\begin{array}{c}
h \\ 
g
\end{array}
\right] (\tau )$ is rational then 
\begin{equation}
Z\left[ 
\begin{array}{c}
h \\ 
g
\end{array}
\right] =Z\left[ 
\begin{array}{c}
h^{d} \\ 
g^{a}
\end{array}
\right]   \label{conjugation property}
\end{equation}
where $(d,pk)=1$ and $(a,p)=1.$
\end{thm}

\smallskip

\textbf{Proof.} Each head character $a_{g,mp-j}(h)$ of (\ref{Z(h,Fricke)})
is rational by assumption. This implies that $%
a_{g,mp-j}(h)=a_{g,mp-j}(h^{d})\,$ for each $d$ such that $(d,pk)=1$ e.g. 
\cite{L}. Hence $Z\left[ 
\begin{array}{c}
h \\ 
g
\end{array}
\right] =Z\left[ 
\begin{array}{c}
h^{d} \\ 
g
\end{array}
\right] \,$since $\phi _{g}(h^{d})=1$. Furthermore, for $p=3,5,7$, $g=p+$ is
conjugate to $g^{a}\,$in ${\bf M}$ for all $(a,p)=1$ so that $g=xg^{a}x^{-1}$
for some $x\in {\bf M}$. Clearly $C_{g}=C_{g^{a}}=xC_{g}x^{-1}\,$ so that $%
xhx^{-1}\in C_{g}$ and $h$ is a member of a $(1,G_{g^{a}})\,$ class of $%
C_{g^{a}}$. Hence Theorem \ref{theorem_conjugation} implies (\ref{general
conjugation}) and $\phi _{g}(xhx^{-1})=1$. Using (\ref{Z2}) we then find $%
Z\left[ 
\begin{array}{c}
h \\ 
g^{a}
\end{array}
\right] =Z\left[ 
\begin{array}{c}
h \\ 
g
\end{array}
\right] $ since $\phi _{g^{a}}(h)=1$. $\qed $

\begin{cor}
\label{corollary_gamma0} Consider $g=p+\ $for $p=2,3,5,7$ with $o(h)=pk$ for 
$k=1$ or $k$ prime. If $Z\left[ 
\begin{array}{c}
h \\ 
g
\end{array}
\right] (\tau )$\ is rational then $\Gamma _{0}^{0}(pk,p)\subseteq \Gamma
_{h,g}$\ i.e. $\Gamma _{0}(p^{2}k)\subseteq \tilde{\Gamma}_{h,g}$.$\,$\ 
\end{cor}

\textbf{Proof.} For $\gamma =\left( 
\begin{array}{ll}
a & b \\ 
c & d
\end{array}
\right) \in \Gamma _{0}^{0}(pk,p)$ we have $b=0\ {\rm mod}\ p\,$and $c=0\ 
{\rm mod}\ pk\,$ as described in Appendix A with $ad-bc=1$. Hence the RHS
of (\ref{SLT}) gives $Z\left[ 
\begin{array}{c}
h^{d} \\ 
g^{a}
\end{array}
\right] (\tau )$ with $(d,pk)=1$ and $(a,p)=1$. The result then follows from
Theorem \ref{theorem_GMF_power_conjugation}. $\qed $

\textbf{Note.} We conjecture that a more general version of Theorem \ref
{theorem_GMF_power_conjugation} holds, namely, if a GMF is rational then (%
\ref{conjugation property}) holds with $(d,o(h))=1$ and $(a,o(g))=1$ so that
the GMF is $\Gamma _{0}^{0}(o(h),o(g))$ invariant.

Theorem \ref{theorem_GMF_power_conjugation} further implies that the
singularity structure is restricted to the Fricke classes of $\langle
g,h\rangle \,$ inequivalent under the equivalence relation $g^{A}h^{B}\sim $ 
$g^{aA}h^{dB}\,$ where $(d,pk)=(a,p)=1$ for $(A,B)=1$. We call the
equivalence classes under the relation $g^{A}h^{B}\sim $ $g^{aA}h^{dB}\,$
for $(d,pk)=1$ and $(a,p)=1$ the class structure of the commuting pair $g,h$%
. We will see below that in all cases considered the class structure
together with Theorems \ref{theorem_proj1}, \ref{theorem_proj2} and \ref
{theorem_GMF_power_conjugation} uniquely determine the GMF and its genus
zero property.

\section{The Genus Zero Property for Some Rational GMFs for $g=p+$}

We now come to the main purpose of this paper which is to demonstrate the
genus zero property for Generalised Moonshine Functions (GMFs) (\ref{Z}) for
rational GMFs with $g=p+$ and $o(h)=pk$ for $k=1$ or $k$ prime. We will show
how the orbifold considerations discussed in the previous sections allow us
to demonstrate that either the given Fricke singularity structure is
inconsistent or else all singularities of the GMF can be identified under
some genus zero fixing group for which the GMF is a hauptmodul. We
conjecture that the method described can be used to demonstrate the genus
zero property for general GMFs including irrational cases. We begin by
discussing four examples of possible Fricke class structures for $\langle
g,h\rangle $ where in the first two examples we demonstrate the genus zero
property and in the second two examples, demonstrate that the given Fricke
class structure is impossible. These examples of class structures are by no
means exhaustive but are nevertheless of general applicability.

\begin{thm}
\label{theorem_h_F_NF} For $g=p+$ for $p=2,3,5,7$ with $o(h)=pk$ for $k=1$
or $k$ prime with $\phi _{g}(h)=1$ then the following class structures give
rise to a rational GMF with genus zero fixing group $\tilde{\Gamma}_{h,g}$
for $Z\left[ 
\begin{array}{c}
h \\ 
g
\end{array}
\right] (p\tau )$ as follows:

\begin{description}
\item[I]  $g$ Fricke and all other classes Non-Fricke then $\tilde{\Gamma}%
_{h,g}=p^{2}k-$.

\item[II]  $g,h$ Fricke and all other classes Non-Fricke then $\tilde{\Gamma}%
_{h,g}=p^{2}k+p^{2}k$.
\end{description}

\textbf{Note.} See Appendix A for the modular group notation.
\end{thm}

\textbf{Proof.} \textbf{I.} Here $Z\left[ 
\begin{array}{c}
h \\ 
g
\end{array}
\right] (p\tau )$ has a unique pole at the cusp $\tau =i\infty $ ($q=0$)
which is invariant under $\Gamma _{0}(p^{2}k)$, from Corollary \ref
{corollary_gamma0}, and is therefore the hauptmodul for $\Gamma _{0}(p^{2}k)$
which must be of genus zero. This restricts $p^{2}k$ to the possible values $%
4,8,9,12,18,25$.

\textbf{\ II.} Here we have singular behaviour in both the $h$ and $g$
twisted sectors, corresponding to cusps at $\tau =0$ and $\tau =i\infty $.
From Theorem \ref{theorem_proj1} (i) using $\phi _{g}(h)=1$ it follows that $%
\mathcal{V}_{\text{orb}}^{\langle g,h\rangle }\simeq \mathcal{V}^{\natural }$%
. Using Theorem \ref{theorem_proj1} (iii) with $u=g$, $v=h$ it follows that
there exists a unique $A$ such that $g^{A}h$ is Fricke of order $pk$ where $%
\phi _{g^{A}h}(g)=1$. But $g^{A}h$ must be therefore conjugate to $h\,$ i.e. 
$A=0\ {\rm mod}\ p$ and $\phi _{h}(g)=1$. Therefore 
\begin{equation}
Z\left[ 
\begin{array}{c}
g^{-1} \\ 
h
\end{array}
\right] (\tau )=q^{-1/pk}+0+O(q^{1/pk})\text{.}
\end{equation}
Define $\hat{W}_{k}$: $\tau \rightarrow -1/k\tau $ which is conjugate to the
Fricke involution $W_{p^{2}k}$ for $\Gamma _{0}(p^{2}k)$ (see Appendix A). $%
\hat{W}_{k}$ normalizes $\Gamma _{0}^{0}(pk,p)$ and interchanges the cusps
at $0$ and $i\infty $.

Let $f(\tau )=Z\left[ 
\begin{array}{c}
h \\ 
g
\end{array}
\right] (\hat{W}_{k}\tau )-Z\left[ 
\begin{array}{c}
h \\ 
g
\end{array}
\right] (\tau )=\ Z\left[ 
\begin{array}{c}
g^{-1} \\ 
h
\end{array}
\right] (k\tau )-Z\left[ 
\begin{array}{c}
h \\ 
g
\end{array}
\right] (\tau )=\ O(q^{1/p})$. $f(\tau )$ is $\Gamma _{0}^{0}(pk,p)$
invariant without poles on ${\bf H}\,\,/\Gamma _{0}^{0}(pk,p)$ and hence is
constant and equal to zero. Hence $\Gamma _{h,g}=\langle \Gamma
_{0}^{0}(pk,p),\hat{W}_{k}\rangle $ i.e. $\tilde{\Gamma}_{h,g}=p^{2}k+p^{2}k$%
. This restricts $p^{2}k$ to the possible values $4,8,9,12,18,20,25,27$. $%
\qed $

\begin{thm}
\label{theorem_impossible} For $g=p+$ for $p=2,3,5,7$ and $o(h)=pk$ for $k=1$
or $k$ prime with $\phi _{g}(h)=1$ the following class structures are
impossible for $k>1.$

\begin{description}
\item[I]  \textbf{\ }$g\,$, $gh\,$ and $h$ Fricke where $g$ is the only
Fricke class of order $p$.

\item[II]  $g$, $gh^{k}$ and $h$ Fricke where $h$ is the only Fricke class
of order $pk$.
\end{description}
\end{thm}

\textbf{Proof. I}\textit{\ }From Theorem \ref{theorem_proj1} (i) $\mathcal{V}%
_{\text{orb}}^{\langle g,h\rangle }\simeq \mathcal{V}^{\natural }$ and from
(iii) with $u=gh$ Fricke and $v=h^{k}$ of order $p$ it follows that $%
u^{A}v=(gh)^{A}h^{k}$ is Fricke for some $A$ with $o(u^{A}v)=p$. Hence $%
u^{A}v=g^{A}h^{A+k}$ must be conjugate to $g$ being the only Fricke class of
order $p$ i.e. $A=-k\ {\rm mod}\ pk\,$ so that $u^{A}v=g^{-k}$ with $\phi
_{g^{-k}}(gh)=1$. However from (\ref{FA}) we also have $\phi
_{g^{-k}}(g^{-k})=\omega _{p}$ which implies that $(\phi
_{g^{-k}}(h))^{k}=\phi _{g^{-k}}((gh)^{k})\phi _{g^{-k}}(g^{-k})=\omega _{p}$%
. But Theorem \ref{theorem_GMF_power_conjugation} implies that since $\phi
_{g}(h)=1$ then $\phi _{g^{-k}}(h)=1\,$ leading to a contradiction. Hence
this conjugacy class structure is impossible.

\textbf{II} Again $\mathcal{V}_{\text{orb}}^{\langle g,h\rangle }\simeq 
\mathcal{V}^{\natural }$. As in Theorem \ref{theorem_h_F_NF} II we can
firstly show that $\phi _{h}(g)=1$. Then applying Theorem \ref{theorem_proj1}
(iii) with $u=gh^{k}$ Fricke and $v=h$ it follows that $u^{A}v=(gh^{k})^{A}h$
is Fricke for some $A$ where $o(u^{A}v)=pk$. Hence $u^{A}v=g^{A}h^{kA+1}$
must be conjugate to $h$ being the only Fricke class of order $pk$ i.e. $%
A=0\ {\rm mod}\ p$ and $u^{A}v=h$ with $\phi _{h}(gh^{k})=1$. But from (%
\ref{FA}) we have $\phi _{h}(h)=\omega _{pk}$ so that $\phi _{h}(g)=(\phi
_{h}(h))^{-k}=\omega _{p}^{-1}$. This contradicts our earlier statement that 
$\phi _{h}(g)=1$. Hence this class structure is also impossible. $\qed $

\subsection{Case A: Rational GMFs for $o(h)=p$, $k=1$}

The analysis when $o(g)=o(h)=p$ with prime $p$ has been previously reported
in \cite{T3}. The class structure of $\langle g,h\rangle $ is now determined
by $g$, $h$ and $gh$. The cusps of ${\bf H}/\Gamma _{0}^{0}(p)$ are at $%
\{i\infty ,0,1,2,...,p-1\}$ corresponding to the sectors twisted by $g$, $h$
and $g^{s}h\sim gh$ ($s=1,2,...,p-1$) respectively with $\Gamma
_{0}^{0}(p)\equiv \Gamma _{0}^{0}(p,p)$ \cite{T1}, \cite{T3}. There are then
four possible cases that may occur.

\begin{thm}
\label{theorem_h_p} For $g=p+$ and $h$ of order $p$ with $\phi _{g}(h)=1$
each possible class structure gives rise to a genus zero fixing group $%
\tilde{\Gamma}_{h,g}$ for rational $Z\left[ 
\begin{array}{c}
h \\ 
g
\end{array}
\right] (p\tau )$ as follows:

\begin{description}
\item[I.A]  $g$ Fricke. $\tilde{\Gamma}_{h,g}=p^{2}-$.

\item[II.A]  $g$, $h$ Fricke. $\tilde{\Gamma}_{h,g}=p^{2}+$.

\item[III.A]  $g$, $gh$ Fricke. $\tilde{\Gamma}_{h,g}=p-$.

\item[IV.A]  $g$, $h$ and $gh$ Fricke. $\tilde{\Gamma}_{h,g}=p|p-$ for $p=2,3
$ or $5||5+$ for $p=5$.
\end{description}

\textbf{Note.} Here every element of $\langle g,h\rangle $ is Non-Fricke
unless otherwise stated. Table 3 below shows all possible examples and the
GMFs actually observed. See Appendix A for the modular group notation.
\end{thm}

\smallskip \textbf{Proof. I.A} This follows from Theorem \ref{theorem_h_F_NF}
I. This restricts the possible values for $p$ to $p=2$, $3$ and $5$. $p=5$
does not arise in practice - see Appendix B \cite{Q1},\cite{Q2}.

\textbf{II.A} This follows from Theorem \ref{theorem_h_F_NF} II. This
restricts the possible values for $p$ to $p=2$, $3$, $5$, $7$. $p=7$ does
not arise in practice - see Appendix B \cite{Q1},\cite{Q2}.

\textbf{III.A }Here we have singular behaviour in the $g$ and $gh$ twisted
sectors. Consider the $SL(2,{\bf Z})$ transformation $\gamma _{p}\equiv
STS=\left( 
\begin{array}{cc}
1 & 0 \\ 
-1 & 1
\end{array}
\right) $ which normalises $\Gamma (p)\equiv \Gamma (p,p)$ and is of order $%
p $ in $\Gamma (p)$. We show that the cusps $\{i\infty ,1,2,...,p-1$\} are
identified under $\gamma _{p}$. From (\ref{SLT}) we have 
\begin{equation}
Z\left[ 
\begin{array}{c}
h \\ 
g
\end{array}
\right] (\gamma _{p}\tau )=\phi _{gh}(h)q^{-1/p}+0+O(q^{1/p})\text{.}
\label{E2}
\end{equation}
First we prove that $\phi _{gh}(h)=1$. According to Theorem \ref
{theorem_proj2} (i) $\mathcal{V}_{\text{orb}}^{\langle g,h\rangle }\simeq 
\mathcal{V}^{\Lambda }$. From Theorem \ref{theorem_proj2} (ii) with $u=g$
Fricke and $v=h$ it follows that $\phi _{gh}(g)\neq 1$ so that $\phi
_{gh}(g)=\omega _{p}^{s}\,$for some $s\neq 0\ \mathrm{mod}\ p$. However with 
$u=g^{a}h$ Fricke for $a\neq 0\ \mathrm{mod}\ p$ and $v=h$ it also follows
that $\phi _{gh}(g^{a}h)\neq 1$. (If $a=1\ \mathrm{mod}\ p$ then from (\ref
{FA}) $\phi _{gh}(gh)=\omega _{p}$ whereas if $a\neq 1\ \mathrm{mod}\ p$
then we can choose $A$ such that $aA=1\ \mathrm{mod}\ p$ and $B=(a-1)A\ 
\mathrm{mod}\ p$, $B\neq 0\ \mathrm{mod}\ p$, so that $%
(g^{a}h)^{A}(h)^{B}=gh $ and then apply Theorem \ref{theorem_proj2} (ii)).
Let $\phi _{gh}(h)=\omega _{p}^{r}$ then $\phi _{gh}(g^{a}h)=\omega
_{p}^{as+r}\neq 1.$ If $r\neq 0\ \mathrm{mod}\ p$ then since $s\neq 0\ 
\mathrm{mod}\ p$ we can always find $a$ such that $as=-r\ \mathrm{mod}\ p$
which is a contradiction. Hence $\phi _{gh}(h)=1$.

From (\ref{E2}) it follows that $f(\tau )\equiv Z\left[ 
\begin{array}{c}
h \\ 
g
\end{array}
\right] (\gamma _{p}\tau )-Z\left[ 
\begin{array}{c}
h \\ 
g
\end{array}
\right] (\tau )=\ O(q^{1/p})$ is $\Gamma (p)$ invariant and is non-singular
at $q=0$. One can easily check that the other cusps of $f$ are similarly
non-singular so that $f$ is holomorphic on ${\bf H}\,\,/\Gamma (p)$ and is
therefore zero. Hence $\Gamma _{h,g}=\langle \Gamma _{0}^{0}(p),\gamma
_{p}\rangle =\Gamma _{0}^{0}(1,p)$ invariant i.e. $\tilde{\Gamma}%
_{h,g}=\Gamma _{0}(p)$. This is possible for $p=2,3,5,7$. Note that for $p=7 
$ this is the only rational class.

\textbf{IV.A }In this case all twisted sectors are Fricke and all cusps are
singular. For $p=2$, $3$ and $5$ there is exactly one rational class of
order $p$ in $G_{p+}$ not already identified in cases I.A to III.A above. We
consider each $p$ in turn.

Consider $p=2$. From Theorem \ref{theorem_proj1} (i) $\mathcal{V}_{\text{orb}%
}^{\langle g,h\rangle }\simeq \mathcal{V}^{\natural }$ and from (iii) with $%
u=h$ and $v=g$ there exists a unique $A\ \mathrm{mod}\ 2$ such that $h^{A}g$
is Fricke and where $\phi _{h^{A}g}(h)=1$ i.e. $A=0$ and $\phi _{gh}(h)=-1$.
From (\ref{FA}) $\phi _{gh}(gh)=-1$ and so $\phi _{gh}(g)=1$. Similarly,
from Theorem \ref{theorem_proj1} (iii) with $u=g$ and $v=h$ there exists a
unique $A\ \mathrm{mod}\ 2$ such that $g^{A}h$ is Fricke and $\phi
_{g^{A}h}(g)=1$ i.e. $A=1$ and $\phi _{h}(g)=-1$. From (\ref{FA}) $\phi
_{h}(h)=-1$ and so $\phi _{h}(gh)=-1$. We summarise these phases in Table 2.

\smallskip

\begin{tabular}{|l|l|l|l|}
\hline
Phase & $g$ & $h$ & $gh$ \\ \hline
$\phi _{g}$ & $-1$ & $1$ & $-1$ \\ \hline
$\phi _{gh}$ & $1$ & $-1$ & $-1$ \\ \hline
$\phi _{h}$ & $-1$ & $-1$ & $1$ \\ \hline
\end{tabular}

Table 2.

Using (\ref{SLT}) we then find that $Z\left[ 
\begin{array}{c}
h \\ 
g
\end{array}
\right] (\tau )$ is also invariant under $ST$ of order $3$ which normalises $%
\Gamma (2)$ and so is the hauptmodul for $\langle \Gamma
_{0}^{0}(2),ST\rangle $ of level two and index two in $SL(2,{\bf Z})$ i.e. $%
\tilde{\Gamma}_{h,g}=2|2$.

For $p=3,5$ we can show that $\phi _{h}(g)=1$. From Theorem \ref
{theorem_GMF_power_conjugation}, $Z\left[ 
\begin{array}{c}
h \\ 
g
\end{array}
\right] =Z\left[ 
\begin{array}{c}
h \\ 
g^{-1}
\end{array}
\right] $ which implies after an $S$ transformation that $Z\left[ 
\begin{array}{c}
g^{-1} \\ 
h
\end{array}
\right] =Z\left[ 
\begin{array}{c}
g \\ 
h
\end{array}
\right] $ so that $\phi _{h}(g)=(\phi _{h}(g))^{-1}$. But $\phi _{h}(g)\in
\langle \omega _{p}\rangle \ $implies $\phi _{h}(g)=1$ for $p>2$. Hence $%
Z\left[ 
\begin{array}{c}
g \\ 
h
\end{array}
\right] (\tau )=\ q^{-1/p}+0+O(q^{1/p})$. For $p=3$ we can consider Theorem 
\ref{theorem_proj1} (iii) with $u=g^{2}h$ and $v=g^{2}.$ Then $\phi
_{g^{2A+2}h^{A}}(g^{2}h)=1$ for a unique $A$. However $\phi
_{g^{2}}(g^{2}h)=\omega _{3}^{2}$ and $\phi _{h^{2}}(g^{2}h)=\omega _{3}^{2}$
and hence $A=1$ so that $\phi _{gh}(g^{2}h)=1$ and hence $\phi _{gh}(g)=\phi
_{gh}(h)=\omega _{3}^{2}$ since $\phi _{gh}(gh)=\omega _{3}$ from (\ref{FA}%
). Hence the phase residues of all the singular cusps are now known. For $%
p=5 $ we observe that $h$ is an element of the $5A$ class of $G_{g}=\mathrm{%
HN}$ being the only remaining available rational class not identified in
cases I.A to III.A\ above. But $\phi _{h}(g)=1$ implies that $g$ is an
element of the $5A$ class of the isomorphic group $G_{h}=\mathrm{HN}$. Hence 
$Z\left[ 
\begin{array}{c}
g \\ 
h
\end{array}
\right] =Z\left[ 
\begin{array}{c}
h \\ 
g
\end{array}
\right] $ which leads to $\phi _{gh}(g)=\phi _{gh}(h)=\omega _{5}^{3}$. A
similar argument can also be used for $p=3$. Thus we find $\phi
_{gh}(g)=\phi _{gh}(h)=\omega _{p}^{(p+1)/2}$ for $p=3,5$.

For $p=3$, $\gamma _{2}=T^{-1}ST$ is of order 2 and interchanges the cusps $%
\{\infty ,0\}\leftrightarrow \{2,1\}$ whereas $S$ interchanges $\{\infty
,1\}\leftrightarrow \{0,2\}$. We then find using the given phase residues
that $\Gamma _{h,g}=\langle \Gamma _{0}^{0}(3),\gamma _{2},S\rangle $ which
is of level $3$ and index $3$ in $SL(2,{\bf Z})$ i.e. $\tilde{\Gamma}%
_{h,g}=3|3$.

For $p=5$, let $\gamma _{3}=TST^{3}$ which is of order $3$ in $\Gamma
_{0}^{0}(5)$ and cyclically permutes the cusps $\{\infty ,0\}\rightarrow
\{1,4\}\rightarrow \{2,3\}$ whereas $S$ interchanges the cusps $\{\infty
,1,2\}\leftrightarrow \{0,4,3\}$. Hence $\Gamma _{h,g}=\langle \Gamma
_{0}^{0}(5),S,\gamma _{3}\rangle $ which is of level $5$ and index $5$ in $%
SL(2,{\bf Z})$ i.e. $\tilde{\Gamma}_{h,g}=$ $5||5+$. $\qed $

We summarise cases I--IV.A in Table 3 where all possible examples are
indicated. Examples that do not arise in explicit calculations are marked
with an asterix and all Fricke classes are in boldface.\smallskip

\begin{tabular}{|l|l|l|l|l|}
\hline
Case & $h$ & $gh$ & $\tilde{\Gamma}_{h,g}$ & Examples \\ \hline
I.A & $p-$ & $p-$ & $\Gamma _{0}(p^{2})$ & $p=2,3,5^{*}$ \\ \hline
II.A & $\mathbf{p+}$ & $p-$ & $\Gamma _{0}(p^{2})+$ & $p=2,3,5,7^{*}$ \\ 
\hline
III.A & $p-$ & $\mathbf{p+}$ & $\Gamma _{0}(p)-$ & $p=2,3,5,7$ \\ \hline
IV.A & $\mathbf{p+}$ & $\mathbf{p+}$ & $p|p-$ or $p||p+$ & $2|2-$, $3|3-$, $%
5||5+$ \\ \hline
\end{tabular}

Table 3. * indicates that no such GMF occurs. Fricke classes are in boldface

\subsection{Case B: Rational GMFs for $o(h)=p^{2}$, $k=p$}

The class structure of $\langle g,h\rangle $ is now described by $g$, $h$, $%
gh$ and $gh^{p}$ where from Theorem \ref{theorem_gamma0_cusps} of Appendix A
the cusps of ${\bf H}/\Gamma _{0}^{0}(p^{2},p)$ are at $\{i\infty $, $0$, $s$%
, $\frac{s}{p}\}$ for $s=1,2,...,p-1$ corresponding to the sectors twisted
by $g$, $h$, $g^{s}h\sim gh$ and $g^{s}h^{p}\sim gh^{p}$ for $s=1$,$2$,...,$%
p-1$ respectively. Furthermore no such $h\in C_{g}$ exists for $p=7$ whereas
for $p=5$ only two algebraically conjugate classes of order $25$ exist with
irrational characters and irrational GMFs. Hence we consider $p=2,3$ only.
There are then eight possible cases that may occur as follows.

\begin{thm}
\label{theorem_h_p^2} For $p=2,3$ with $g=p+$ and $h$ of order $p^{2}$ with $%
\phi _{g}(h)=1$ some possible class structures gives rise to a genus zero
fixing group $\tilde{\Gamma}_{h,g}$ for rational $Z\left[ 
\begin{array}{c}
h \\ 
g
\end{array}
\right] (p\tau )$ while others are impossible as follows:

\begin{description}
\item[I.B]  $g$ Fricke. $\tilde{\Gamma}_{h,g}=p^{3}-$.

\item[II.B]  $g$, $h$ Fricke. $\tilde{\Gamma}_{h,g}=p^{3}+$.

\item[III.B]  $g$, $gh^{p}$ Fricke. $\tilde{\Gamma}_{h,g}=p^{2}-$ or $p^{2}|p
$.

\item[IV.B]  $g$, $gh$ and $gh^{p}$ Fricke. Impossible.

\item[V.B]  $g$, $h$ and $gh$ Fricke. Impossible.

\item[VI.B]  $g$, $h$ and $gh^{p}$ Fricke. Impossible.

\item[VII.B]  $g$ and $gh$ Fricke. $\tilde{\Gamma}_{h,g}=8\frac{1}{2}+$ for $%
p=2$. Impossible for $p=3$.

\item[VIII.B]  $g$, $h$, $gh$ and $gh^{p}$ Fricke. $\tilde{\Gamma}_{h,g}=4|2+
$ or $4|2+2^{\prime }$ for $p=2$, $\tilde{\Gamma}_{h,g}=$ $9|3+$ or $9||3+$
for $p=3$.
\end{description}

\textbf{Note.} Every element of $\langle g,h\rangle $ is Non-Fricke unless
otherwise stated. Table 5 below shows all possible examples and the GMFs
actually observed. See Appendix A for the modular group notation.
\end{thm}

\textbf{Proof. I.B} This follows from Theorem \ref{theorem_h_F_NF} I.

\textbf{II.B} This follows from Theorem \ref{theorem_h_F_NF} II.

\textbf{III.B }We have singular behaviour in the $g$ and $gh^{p}$ twisted
sectors. Let $\phi _{gh^{p}}(h)=\omega _{p}^{b}$ for some $b$ using (\ref{FA}%
). Consider $\gamma _{p}^{(b)}\equiv ST^{p}ST^{b}=\left( 
\begin{array}{cc}
1 & b \\ 
-p & 1-pb
\end{array}
\right) $ which normalises and is of order $p$ in $\Gamma (p^{2},p)$. We
have from (\ref{SLT}) that 
\begin{equation}
Z\left[ 
\begin{array}{c}
h \\ 
g
\end{array}
\right] (\gamma _{p}^{(b)}\tau )=\phi
_{gh^{p}}(h^{1-pb}g^{-b})q^{-1/p}+0+O(q^{1/p})  \label{E6}
\end{equation}
Then $\phi _{gh^{p}}(h^{1-pb}g^{-b})=(\phi _{gh^{p}}(gh^{p}))^{-b}\phi
_{gh^{p}}(h)=1$ due to (\ref{FA}). From (\ref{E6}) it follows that $Z\left[ 
\begin{array}{c}
h \\ 
g
\end{array}
\right] (\gamma _{p}\tau )-Z\left[ 
\begin{array}{c}
h \\ 
g
\end{array}
\right] (\tau )=O(q^{1/p})$ is holomorphic on ${\bf H}\,/\Gamma (p^{2},p)$
and is therefore zero. Thus $\Gamma _{h,g}=\langle \Gamma
_{0}^{0}(p^{2},p),\gamma _{p}^{(b)}\rangle $. If $b=0$ mod $p$ then 
$\Gamma _{h,g}=\Gamma _{0}^{0}(p,p)\ $and $\tilde{\Gamma}_{h,g}=\Gamma
_{0}(p^{2})$ invariant of genus zero. If $b\neq 0$ mod $p$ then $%
\Gamma _{h,g}$ has index $p$ in $\Gamma _{0}^{0}(p,1)=\Gamma _{0}(p)$ and we
find $\tilde{\Gamma}_{h,g}=p^{2}|p-$ of genus zero.

\textbf{IV.B} From (\ref{FA}) $\phi _{gh}(gh)=\omega _{p^{2}}\,$so that $%
\phi _{gh}(h^{p})=$ $\omega _{p}$. From Theorem \ref{theorem_proj2} (ii)
with $u=g$ Fricke and $v=h$ it follows that $\phi _{gh}(g)\neq 1$ so that $%
\phi _{gh}(g)=\omega _{p}^{s}$ for some $s\neq 0\ {\rm mod}\ p$ since $g$
is order $p$. Taking $u=g^{a}h^{p}$ Fricke for any $a\neq 0\ {\rm mod}\ p$
and $v=h$ it also follows that $\phi _{gh}(g^{a}h^{p})\neq 1$ (similarly to
Theorem \ref{theorem_h_p} III.A). Choose $a$ such that $as=-1\ {\rm mod}\ p$
so that $\phi _{gh}(g^{a}h^{p})=\omega _{p}^{as+1}=1$ which is a
contradiction. Hence this class structure is impossible.

\textbf{V.B} From Theorem \ref{theorem_impossible} I this class structure is
impossible.

\textbf{VI.B }From Theorem \ref{theorem_impossible} II this class structure
is impossible.

\textbf{VII.B} From (\ref{FA}) $\phi _{gh}(gh)=\omega _{p^{2}}$. From
Theorem \ref{theorem_proj2} (ii) with $u=g$ Fricke and $v=h$ it follows that 
$\phi _{gh}(g)\neq 1$. For $p=2$ these conditions imply $\phi _{gh}(g)=-1$
and $\phi _{gh}(h)=\omega _{4}^{-1}$. Consider $\gamma _{2}^{(2)}=\left( 
\begin{array}{cc}
1 & -1 \\ 
-1 & 2
\end{array}
\right) \left( 
\begin{array}{cc}
2 & 0 \\ 
0 & 1
\end{array}
\right) $ which normalizes $\Gamma _{0}^{0}(4,2)$ and is of order $2$. We
then find that $Z\left[ 
\begin{array}{c}
h \\ 
g
\end{array}
\right] (\gamma _{2}^{(2)}(\tau ))=Z\left[ 
\begin{array}{c}
gh^{2} \\ 
gh
\end{array}
\right] (2\tau )=q^{-1/2}+O(q^{1/2})$. Hence we find that $\Gamma _{h,g}=$ $%
\langle \Gamma _{0}^{0}(4,2),\gamma _{2}^{(2)}\rangle $ i.e. $\tilde{\Gamma}%
_{h,g}=8\frac{1}{2}+\,$of genus zero.

For $p=3$ we will show that the orbifolding procedure is not consistent.
Since $h^{3}=$ $3-$ then from (\ref{Tastar}) we have $Z\left[ 
\begin{array}{c}
h^{3} \\ 
1
\end{array}
\right] (\tau )=\left[ \frac{\eta (\tau )}{\eta (3\tau )}\right] ^{12}+12$
i.e. $\rho _{h^{3}}^{0}(h^{3})=\mathbf{1}\mathsf{\ }$and $\chi
_{h^{3}}^{0}(h^{3})=12$ from (\ref{Unit}). Since $gh=9-$ then $\mathcal{V}_{%
\text{orb}}^{\langle gh\rangle }\simeq \mathcal{V}^{\natural }$ and the
constant term of $\mathrm{Tr}_{\mathcal{H}_{(gh)^{A}}}(\mathcal{P}%
_{gh}q^{L_{0}-1})\,$ vanishes for all $A$. Taking $A=3$ we have $\frac{1}{9}%
\sum_{A=0}^{8}\chi _{h^{3}}^{0}((gh)^{A})=\ \frac{1}{3}[12+\chi
_{h^{3}}^{0}(gh)+\chi _{h^{3}}^{0}((gh)^{-1}))]=0$. Hence $\chi
_{h^{3}}^{0}(gh)+\chi _{h^{3}}^{0}((gh)^{-1})=-12$. Since $h=9-$ then $%
\mathcal{V}_{\text{orb}}^{\langle h\rangle }\simeq \mathcal{V}^{\Lambda }$
and the constant term of $\mathrm{Tr}_{\mathcal{H}_{\text{orb}}^{\langle
h\rangle }}(q^{L_{0}-1})\,$is $24$. Hence $\frac{1}{9}\sum_{A=0}^{8}%
\sum_{B=1}^{8}\chi _{h^{B}}^{0}(h^{A})=24$. But $\rho _{h}^{0}(h)=\mathbf{1\,%
}$with $\chi _{h}^{0}(h)=3$ from (\ref{Unit}). From this it follows that $%
\chi _{h^{3}}^{0}(h)+\chi _{h^{3}}^{0}(h^{2})=-3$.

From Theorem \ref{theorem_h_p} I.A we have $Z\left[ 
\begin{array}{c}
h^{3} \\ 
g
\end{array}
\right] (3\tau )=\left[ \frac{\eta (\tau )}{\eta (9\tau )}\right] ^{3}+3$
and after an $S$ transformation we have $Z\left[ 
\begin{array}{c}
g^{2} \\ 
h^{3}
\end{array}
\right] (\tau /3)=3+27\left[ \frac{\eta (\tau )}{\eta (\tau /9)}\right]
^{3}=3+O(q^{1/9})$ i.e. $\chi _{h^{3}}^{0}(g)=\chi _{h^{3}}^{0}(g^{2})=3$.
Using this information we find an inconsistency in the orbifolding with
respect to $\langle g,h\rangle $ where in particular, $\chi _{h^{3}}^{0}(%
\mathcal{P}_{\langle g,h\rangle })=\sum_{u\in \langle g,h\rangle }\chi
_{h^{3}}^{0}(u)=-1$ which must be a non-negative quantity being the
dimension of the $h^{3}\,$twisted vacuum space invariant under $\langle
g,h\rangle $. Hence this class structure is impossible for $p=3$.

\textbf{VIII.B} We firstly consider $p=2$. From Theorem \ref{theorem_proj1}
(i) it follows that $\mathcal{V}_{\text{orb}}^{\langle g,h\rangle }\simeq 
\mathcal{V}^{\natural }$ and from (ii) with $u=gh$ Fricke and $v=g$ it
follows that there exists a unique $A$ such that $u^{A}v=g^{A+1}h^{A}=2+$
with $\phi _{g^{A+1}h^{A}}(gh)=1$. This implies that $g^{A+1}h^{A}$ is
conjugate to either $g$ or $gh^{2}$ i.e. $A=0$ or $2\ {\rm mod}\ 4$. If $%
A=0\ {\rm mod}\ 4$ then $\phi _{g}(gh)=1$ which contradicts the known
values $\phi _{g}(h)=1$ and $\phi _{g}(g)=-1$. Hence $A=2\ {\rm mod}\ 4$
and $\phi _{gh^{2}}(gh)=1$, $\phi _{gh^{2}}(gh^{2})=-1$ so that $\phi
_{gh^{2}}(g)=\phi _{gh^{2}}(h)=-1$.

Similarly from Theorem \ref{theorem_proj1} (iii) with $u=gh^{2}$ Fricke and $%
v=h^{-1}$ it follows that there exists a unique $A$ such that $%
u^{A}v=g^{A}h^{2A-1}=4+$ with $\phi _{g^{A}h^{2A-1}}(gh^{2})=1$. This
implies that $g^{A}h^{2A-1}$ must be conjugate to either $h$ or $gh$ i.e. $%
A=0$ or $1\ {\rm mod}\ 2$. Suppose $A=0$ ${\rm mod}\ 2$. Then $\phi
_{h}(gh^{2})=1$ and so $\phi _{h}(h)=i$ $\,$implies that $\phi _{h}(g)=-1$.
Furthermore, from Theorem \ref{theorem_proj1} (iii) we have $\phi
_{gh}(gh^{2})\neq 1$ which means $\phi _{gh}(h)\neq -i$ since $\phi
_{gh}(gh)=i$. Hence $\phi _{gh}(h)=i$ (since $gh=4+$) and hence $\phi
_{gh}(g)=1$. If on the other hand $A=1$ ${\rm mod}\ 2$ then $\phi
_{gh}(gh^{2})=1$ and so $\phi _{gh}(gh)=i$ $\,$implies that $\phi
_{gh}(h)=-i $ and $\phi _{gh}(g)=-1$. Furthermore, from Theorem \ref
{theorem_proj1} (iii) we have $\phi _{h}(gh^{2})=-1$ (since $gh^{2}$ is of
order $2$) which means $\phi _{h}(g)=1$ since $\phi _{h}(h)=i$. Hence all
residues of singular cusps are determined for $p=2$ for both values of $A.\,$
In particular note that $\phi _{h}(g)=-1$ for $A=0$ ${\rm mod}\ 2$ and $%
\phi _{h}(g)=1$ for $A=1$ ${\rm mod}\ 2$.

Now we can apply an analysis similar to III.B and using the above
information to obtain $\langle \Gamma _{0}^{0}(4,2),\gamma
_{2}^{(-1)}\rangle \ $invariance with $b=-1$, $p=2$ in $\gamma _{p}^{(b)}$
of III.B. We then find that $\hat{w}_{2}:\tau \rightarrow -1/2\tau $ either
fixes or negates $Z\left[ 
\begin{array}{c}
h \\ 
g
\end{array}
\right] (\tau )$ depending on the sign of $\phi _{h}(g)$. Hence $\Gamma
_{h,g}=\langle \Gamma _{0}^{0}(4,2),\gamma _{2}^{(-1)},\hat{w}_{2}\rangle $
or $\langle \Gamma _{0}^{0}(4,2),\gamma _{2}^{(-1)},\hat{w}_{2}^{\prime
}\rangle $ of genus zero i.e. $\tilde{\Gamma}_{h,g}=4|2+$ or $4|2+2^{\prime
} $ where $2^{\prime }$ indicates that the GMF is negated by $w_{2}$.

Consider $p=3$. As in IV.A we find that the rationality of the GMF implies
that $\phi _{h}(g)=1$. Next $\mathcal{V}_{\text{orb}}^{\langle g,h\rangle
}\simeq \mathcal{V}^{\natural }$ and from Theorem \ref{theorem_proj1} (ii)
with $u=gh$ Fricke and $v=g$ it follows that there exists a unique $A$ such
that $u^{A}v=g^{A+1}h^{A}$ is Fricke with $\phi _{g^{A+1}h^{A}}(gh)=1$. This
implies that $g^{A+1}h^{A}$ is conjugate to either $g$ or $gh^{3}$ i.e. $%
A=0\ {\rm mod}\ 3$. We let $b=-A/3$. If $b=0\ {\rm mod}\ 3$ then $\phi
_{g}(gh)=1$ which contradicts the known values $\phi _{g}(h)=1$ and $\phi
_{g}(g)=\omega _{3}$. If $b=1\ {\rm mod}\ 3$ then $\phi _{gh^{6}}(gh)=1$
and $\phi _{gh^{6}}(gh^{6})=\omega _{3}$ implies that $\phi
_{gh^{6}}(g)=\phi _{gh^{6}}(h^{-1})=\omega _{3}$ so that $\phi
_{gh^{3}}(g)=\phi _{gh^{3}}(h)=\omega _{3}$. If $b=-1\ {\rm mod}\ 3$ then $%
\phi _{gh^{3}}(gh)=1$ and $\phi _{gh^{3}}(gh^{3})=\omega _{3}$ implies that $%
\phi _{gh^{3}}(g)=\omega _{3}$ and $\phi _{gh^{3}}(h)=\omega _{3}^{-1}$.
Thus we find $\phi _{gh^{3}}(g)=\omega _{3}$ and $\phi _{gh^{3}}(h)=\omega
_{3}^{b}\,$ for $b=\pm 1\ {\rm mod}\ 3$.

\smallskip We now show that $\phi _{gh}(g)=\omega _{3}^{b}$ and $\phi
_{gh}(h)=\omega _{9}^{1-3b}$. Firstly we have $\phi _{gh^{3}}(g^{-b}h)=1$
and hence as in IV.A we find that the rationality of the GMF implies that $%
\phi _{g^{-b}h}(gh^{3})=1$. Conjugating as in (\ref{phi_class_fun}) we then
find that $\phi _{gh}(gh^{-3b})=1$. However $\phi _{gh}(gh)=\omega _{9}$
implies that $\phi _{gh}(g)=\omega _{3}^{b}$ and $\phi _{gh}(h)=\omega
_{9}^{1-3b}$ given that $b=\pm 1\ {\rm mod}\ 3$. Hence all residues of
singular cusps are now determined as follows:\smallskip

\begin{tabular}{|l|l|l|}
\hline
Phase & $g$ & $h$ \\ \hline
$\phi _{g}$ & $\omega _{3}$ & $1$ \\ \hline
$\phi _{h}$ & $1$ & $\omega _{9}$ \\ \hline
$\phi _{gh}$ & $\omega _{3}^{b}$ & $\omega _{9}^{1-3b}$ \\ \hline
$\phi _{gh^{3}}$ & $\omega _{3}$ & $\omega _{3}^{b}$ \\ \hline
\end{tabular}

Table 4.

Similarly to III.B and using the Table 4 we then obtain $\langle \Gamma
_{0}^{0}(9,3),\gamma _{3}^{(b)},\hat{w}_{3}\rangle \ $ invariance with $%
b=\pm 1$ and $p=3$ for $\gamma _{p}^{(b)}$ of III.B. This corresponds to $%
\tilde{\Gamma}_{h,g}=9|3+$ for $b=-1\ {\rm mod}\ 3$ and $\tilde{\Gamma}%
_{h,g}=9||3+$ for $b=1\ {\rm mod}\ 3\,$ both of genus zero. $\qed $

We summarise cases I--VIII.B in Table 5 where all possible examples are
indicated.

\begin{tabular}{|l|l|l|l|l|l|}
\hline
Case & $h$ & $gh$ & $gh^{p}$ & $\tilde{\Gamma}_{h,g}$ & Examples \\ \hline
I.B & $p^{2}-$ & $p^{2}-$ & $p-$ & $p^{3}-$ & $p=2$ \\ \hline
II.B & $\mathbf{p}^{2}\mathbf{+}$ & $p^{2}-$ & $p-$ & $p^{3}+$ & $p=2,3$ \\ 
\hline
III.B & $p^{2}-$ & $p^{2}-$ & $\mathbf{p+}$ & $p^{2}-$,$\ p^{2}|p-$ & $p=2,3$
\\ \hline
IV.B & $p^{2}-$ & $\mathbf{p}^{2}\mathbf{+}$ & $\mathbf{p+}$ & Impossible & -
\\ \hline
V.B & $\mathbf{p}^{2}\mathbf{+}$ & $\mathbf{p}^{2}\mathbf{+}$ & $p-$ & 
Impossible & - \\ \hline
VI.B & $\mathbf{p}^{2}\mathbf{+}$ & $p^{2}-$ & $\mathbf{p+}$ & Impossible & -
\\ \hline
VII.B & $p^{2}-$ & $\mathbf{p}^{2}\mathbf{+}$ & $p-$ & $8\frac{1}{2}+$ & $%
p=2 $ \\ \hline
VIII.B & $\mathbf{p}^{2}\mathbf{+}$ & $\mathbf{p}^{2}\mathbf{+}$ & $\mathbf{%
p+}$ & $4|2+$, $4|2+2^{\prime }$, $9|3+,9||3+$ & $p=2,3$ \\ \hline
\end{tabular}

Table 5. Fricke classes are in boldface.

\subsection{Case C: Rational GMFs for $o(h)=pk$, $k$ prime, $k\neq p$}

For $o(h)=pk$, $k$ prime, $k\neq p$, the class structure of $\langle
g,h\rangle $ is described by $g$, $h$, $gh$, $gh^{p}$, $h^{k}$ and $gh^{k}$.
From Theorem \ref{theorem_gamma0_cusps} of Appendix A the cusps of ${\bf H}%
/\Gamma _{0}^{0}(pk,p)$ are at $\{i\infty ,0,s,\frac{1}{p},\frac{p}{k},\frac{%
s}{k}\}$ for $s=1,2,...,p-1$ corresponding to the sectors twisted by $g$, $h$%
, $g^{s}h\sim gh$, $gh^{p}$, $h^{k}$ and $g^{s}h^{k}\sim gh^{k}$ for $s=1$,$%
2 $,...,$p-1$ respectively. $h^{k}$ and $gh^{k}$ are determined by the power
map formula for $h$ and $gh$. There are then thirteen possible cases that
may occur as follows.\smallskip

\begin{tabular}{|l|l|l|l|l|l|}
\hline
Case & $h$ & $gh$ & $gh^{p}$ & $h^{k}$ & $gh^{k}$ \\ \hline
I.C & $pk-$ & $pk-$ & $pk+p$ & $p-$ & $p-$ \\ \hline
II.C & $\mathbf{pk+pk}$ & $pk-$ & $pk+p$ & $p-$ & $p-$ \\ \hline
III.C & $pk+p$ & $pk-$ & $pk+p$ & $\mathbf{p+}$ & $p-$ \\ \hline
IV.C & $pk+p$ & $pk+k$ & $\mathbf{pk+}$ & $p-$ & $p-$ \\ \hline
V.C & $\mathbf{pk+}$ & $pk+k$ & $\mathbf{pk+}$ & $\mathbf{p+}$ & $p-$ \\ 
\hline
VI.C & $pk-$ & $pk+p$ & $pk+p$ & $p-$ & $\mathbf{p+}$ \\ \hline
VII.C & $\mathbf{pk+pk}$ & $pk+p$ & $pk+p$ & $p-$ & $\mathbf{p+}$ \\ \hline
VIII.C & $pk+p$ & $pk+p$ & $pk+p$ & $\mathbf{p+}$ & $\mathbf{p+}$ \\ \hline
IX.C & $pk-$ & $\mathbf{pk+pk}$ & $pk+p$ & $p-$ & $p-$ \\ \hline
X.C & $\mathbf{pk+pk}$ & $\mathbf{pk+pk}$ & $pk+p$ & $p-$ & $p-$ \\ \hline
XI.C & $pk+p$ & $\mathbf{pk+pk}$ & $pk+p$ & $\mathbf{p+}$ & $p-$ \\ \hline
XII.C & $pk+k$ & $\mathbf{pk+}$ & $\mathbf{pk+}$ & $p-$ & $\mathbf{p+}$ \\ 
\hline
XIII.C & $\mathbf{pk+}$ & $\mathbf{pk+}$ & $\mathbf{pk+}$ & $\mathbf{p+}$ & $%
\mathbf{p+}$ \\ \hline
\end{tabular}

Table 6. Possible classes with Fricke classes shown in boldface

\begin{thm}
\label{theorem_h_pk} With $g=p+$ and $h$ of order $pk$ for $k$ prime $k\neq
p\,$ with $\phi _{g}(h)=1$ some class structures give rise to a genus zero
fixing group $\tilde{\Gamma}_{h,g}$ for rational $Z\left[ 
\begin{array}{c}
h \\ 
g
\end{array}
\right] (p\tau )$ while others are impossible as follows:

\begin{description}
\item[I.C]  $g$ Fricke. $\tilde{\Gamma}_{h,g}=p^{2}k-$ .

\item[II.C]  $g$, $h$ Fricke. $\tilde{\Gamma}_{h,g}=p^{2}k+p^{2}k$.

\item[III.C]  $h^{k}$ Fricke. $\tilde{\Gamma}_{h,g}=p^{2}k+p^{2}$.

\item[IV.C]  $gh^{p}$ Fricke. $\tilde{\Gamma}_{h,g}=p^{2}k+k\,$ for $%
p=2,3,5,7$ or $4k+k^{\prime }$ for $p=2$.

\item[V.C]  $g$, $h$, $gh^{p}$\ and $h^{k}$\ Fricke. $\tilde{\Gamma}%
_{h,g}=p^{2}k+$.

\item[VI.C]  $g$\ and $gh^{k}$\ Fricke. $\tilde{\Gamma}_{h,g}=pk-$.

\item[VII.C]  $g$, $h$ and $gh^{k}$\ Fricke. Impossible.

\item[VIII.C]  $g$, $h^{k}$ and $gh^{k}$\ Fricke. $\tilde{\Gamma}_{h,g}=pk|p-
$ for $p=2,3$. Impossible for $p=5,7$.

\item[IX.C]  $g$\ and $gh$\ Fricke. $\tilde{\Gamma}_{h,g}=4k\frac{1}{2}+4k$
when $p=2$. Impossible for $p=3,5,7$.

\item[X.C]  $g$, $h$\ and $gh$\ Fricke. Impossible.

\item[XI.C]  $g$, $gh$\ and $h^{k}$\ Fricke. $\tilde{\Gamma}_{h,g}=\Gamma [9k%
\symbol{126}a]$ when $p=3$. Impossible for $p=2,5,7$.

\item[XII.C]  $g$, $gh$, $gh^{p}$\ and $gh^{k}$\ Fricke. $\tilde{\Gamma}%
_{h,g}=pk+k$.

\item[XIII.C]  All Fricke. $\tilde{\Gamma}_{h,g}=2k|2+k,2k|2+k^{\prime }$
for $p=2$, $\tilde{\Gamma}_{h,g}=pk|p+$ for $p=2,3$ and $\tilde{\Gamma}%
_{h,g}=pk||p+$ for $p=3,5$.
\end{description}

\textbf{Note.} Table 12 below shows all possible examples and the GMFs
actually observed. Every element of $\langle g,h\rangle $ is Non-Fricke
unless otherwise stated. See Appendix A for the modular group notation.
\end{thm}

\textbf{Proof. I.C}\textit{\ }This follows from Theorem \ref{theorem_h_F_NF}
I.

\textbf{II.C} This follows from Theorem \ref{theorem_h_F_NF} II.

\textbf{III.C} Since $h^{k}$\textit{\ }is Fricke and $\phi _{g}(h^{k})=(\phi
_{g}(h))^{k}=1$ $\,$then from Theorem \ref{theorem_proj1} (i) it follows
that $\mathcal{V}_{\text{orb}}^{\langle g,h^{k}\rangle }\simeq \mathcal{V}%
^{\natural }$. Then from (iii) with $u=g$ Fricke and $v=h^{k}$ it follows
that there exists a unique $A$ such that $g^{A}h^{k}$ is Fricke and $\phi
_{g^{A}h^{k}}(g)=1$. But $g^{A}h^{k}$ must be conjugate to $g$ or $h^{k}$
i.e. $A=0\ {\rm mod}\ p$ and $\phi _{h^{k}}(g)=1$. Furthermore since $%
(k,p)=1\,$and $\phi _{h^{k}}(h^{p})\in \langle \omega _{p}\rangle $ it
follows that $\phi _{h^{k}}(h^{p})=1$.

Consider the $SL(2,{\bf Z)}$ transformation $\hat{W}_{p^{2}}=\left( 
\begin{array}{cc}
p & b \\ 
k & pd
\end{array}
\right) $, $p^{2}d-bk=1$, conjugate to $W_{p^{2}}$, the Atkin-Lehner
involution for $\Gamma _{0}(p^{2}k)$ from Appendix A. Then from (\ref{SLT}) $%
Z\left[ 
\begin{array}{c}
h \\ 
g
\end{array}
\right] (\hat{W}_{p^{2}}\tau )-Z\left[ 
\begin{array}{c}
h \\ 
g
\end{array}
\right] (\tau )=\ (\phi _{h^{k}}(h^{-pd}g^{b})-1)q^{-1/p}+0+O(q^{1/p})=\
O(q^{1/p})$ is $\Gamma _{0}^{0}(pk,p)$ invariant without poles on ${\bf H}%
\,\,/\Gamma _{0}^{0}(pk,p)$ and hence is constant and equal to zero.
Therefore $\Gamma _{h,g}=\langle \Gamma _{0}^{0}(pk,p)$, $\hat{W}%
_{p^{2}}\rangle $ is a genus zero group i.e. $\tilde{\Gamma}%
_{h,g}=p^{2}k+p^{2}$.

\textbf{\ IV.C} We have singular behaviour in the $gh^{p}$ and $g$ twisted
sectors corresponding to the cusps $\tau =1/p$ and $\tau =i\infty $.
Consider the transformation $\hat{W}_{k}=\left( 
\begin{array}{cc}
-k & pb \\ 
pk & -kd
\end{array}
\right) =\left( 
\begin{array}{cc}
-1 & pb \\ 
p & -kd
\end{array}
\right) \left( 
\begin{array}{cc}
k & 0 \\ 
0 & 1
\end{array}
\right) $, $\det (\hat{W}_{k})=k$; $kd-bp^{2}=1$. $\hat{W}_{k}$ is conjugate
to $W_{k}=\left( 
\begin{array}{cc}
-k & b \\ 
p^{2}k & -kd
\end{array}
\right) $ the Atkin-Lehner involution for $\Gamma _{0}(p^{2}k)$ of Appendix
A. Then from (\ref{SLT}) $Z\left[ 
\begin{array}{c}
h \\ 
g
\end{array}
\right] (\hat{W}_{k}\tau )=\ \phi
_{gh^{p}}(h^{kd}g^{bp})q^{-1/p}+O(q^{1/p})=\ \phi
_{gh^{p}}(h^{kd})q^{-1/p}+O(q^{1/p})$.

If $p=2$ then $o(h^{k})=2\,\,$and hence $\phi _{gh^{2}}(h^{k})=\pm 1$. Then
analogously to Theorem \ref{theorem_h_p^2} VIII.B we conclude that $Z\left[ 
\begin{array}{c}
h \\ 
g
\end{array}
\right] (\hat{W}_{k}\tau )=\pm Z\left[ 
\begin{array}{c}
h \\ 
g
\end{array}
\right] (\tau )$ and $\tilde{\Gamma}_{h,g}=4k+k$ or $4k+k^{\prime }$ both of
genus zero where $k^{\prime }$ indicates that the GMF is negated by $W_{k}$.
For $p>2$ we will now show that $\phi _{gh^{p}}(h^{k})=1$. From Theorem \ref
{theorem_GMF_power_conjugation} we have $Z\left[ 
\begin{array}{c}
h \\ 
g
\end{array}
\right] (\tau )=Z\left[ 
\begin{array}{c}
h^{1+ak} \\ 
g
\end{array}
\right] (\tau )$ for $(1+ak,pk)=1$. But for all $a$ we have $(1+ak,k)=1$ and
either $(1+k,p)=1$ or $(1+2k,p)=1$ and $1+2k<pk$ i.e. we can choose $a=1$ or 
$2$. Applying $\hat{W}_{k}$ we then find that $\phi _{gh^{p}}(h^{kd})=\phi
_{gh^{p}}(h^{kd(1+ak)})$ i.e. $(\phi _{gh^{p}}(h^{k^{2}}))^{ad}=1$. However
from Theorem \ref{theorem_phi} $\phi _{gh^{p}}(h^{k})\in \langle \omega
_{p}\rangle $ and $(ad,p)=1$ and so $\phi _{gh^{p}}(h^{k})=1$. Hence $%
Z\left[ 
\begin{array}{c}
h \\ 
g
\end{array}
\right] (\hat{W}_{k}\tau )-Z\left[ 
\begin{array}{c}
h \\ 
g
\end{array}
\right] (\tau )=\ O(q^{1/p})$ is without poles on ${\bf H}\,\,/\Gamma
_{0}^{0}(pk,p)$ and hence is constant and equal to zero. So $\Gamma
_{h,g}=\langle \Gamma _{0}^{0}(pk,p)$, $\hat{W}_{k}\rangle $ or $\tilde{%
\Gamma}_{h,g}=p^{2}k+k$ of genus zero for $p=3,5,7$.

\textbf{\ V.C} From Theorem \ref{theorem_phi} (i) $\mathcal{V}_{\text{orb}%
}^{\langle g,h\rangle }\simeq \mathcal{V}^{\natural }$ since $h$ is Fricke.
From (iii) with $u=h^{k}$ Fricke $v=gh^{p}$ $\langle u,v\rangle =\langle
g,h\rangle $ and from (iii) it follows that there exists a unique $A$ such
that $gh^{Ak+p}$ is Fricke with $\phi _{gh^{Ak+p}}(h^{k})=1$. But then $%
gh^{Ak+p}$ must be conjugate to $gh^{p}$. Therefore $A=p$ with $\phi
_{gh^{p}}(h^{k})=1$. As in III.C above and Theorem \ref{theorem_h_p} II.A we
can also show that $\phi _{h^{k}}(h^{p})=1$ and $\phi _{h}(g)=1$. Using the
standard argument we find that $\tilde{\Gamma}_{h,g}=p^{2}k+$ of genus zero.

\textbf{\ VI.C} Consider the $SL(2,{\bf Z})$ transformation $\gamma
_{p}=ST^{k}S$ which is of order $p$ in and normalises $\Gamma (pk,p)$. From (%
\ref{SLT}) we have 
\begin{equation}
Z\left[ 
\begin{array}{c}
h \\ 
g
\end{array}
\right] (\gamma _{p}\tau )=\phi _{gh^{k}}(h)q^{-1/p}+0+O(q^{1/p})
\label{Z(h,g)(gammaptau)}
\end{equation}

First we will prove that $\phi _{gh^{k}}(h)=1$. From Theorem \ref
{theorem_proj2} (i) $\mathcal{V}_{\text{orb}}^{\langle g,h^{k}\rangle
}\simeq \mathcal{V}^{\Lambda }$. With $u=g$ Fricke and $v=h^{k}$ from (ii)
it follows that $\phi _{gh^{k}}(g)=\omega _{p}^{s}\neq 1$. With $%
u=g^{a}h^{k} $ Fricke for $a\neq 0\ {\rm mod}\ p$ and $v=h^{k}$ it also
follows that $\phi _{gh^{k}}(g^{a}h^{k})\neq 1$ (Here as in Theorem \ref
{theorem_h_p} III.A we can choose $A$ and $B$ such that $u^{A}v^{B}=gh^{k}$
for any $a\neq 0\ {\rm mod}\ p$ and apply Theorem \ref{theorem_proj2}
(ii)). Let $\phi _{gh^{k}}(h)=\omega _{p}^{r}$ then $\phi
_{gh^{k}}(g^{a}h)=\omega _{p}^{as+kr}\neq 1.$ If $r\neq 0\ \mathrm{mod}\ p$
then since $s\neq 0\ \mathrm{mod}\ p$ we can always find $a$ such that $%
as=-kr\ \mathrm{mod}\ p$ which is a contradiction. Hence $\phi
_{gh^{k}}(h)=1 $. Now from (\ref{Z(h,g)(gammaptau)}) it follows that $\Gamma
_{h,g}=\langle \Gamma _{0}^{0}(pk,p),\gamma _{p}\rangle \equiv \Gamma
_{0}^{0}(k,p)$ of genus zero and $\tilde{\Gamma}_{h,g}=\Gamma _{0}(pk)$
where all cusps $\tau =i\infty $, $\frac{s}{k}$, ($s=1$,$2$,...,$p-1$) are
identified under $\gamma _{p}$.

\textbf{VII.C} This class structure is impossible as shown in Theorem \ref
{theorem_impossible} II.

\textbf{VIII.C }From Theorem \ref{theorem_h_p} IV.A we must have $p=2,3$ or $%
5$. Let us first consider $p=2$. All Fricke classes are of order $2$ and
since $k\neq p=2$ and using Table 2 of Theorem \ref{theorem_h_p} IV.A we
obtain:\smallskip

\begin{tabular}{|l|l|l|}
\hline
Phase & $g$ & $h$ \\ \hline
$\phi _{g}$ & $-1$ & $1$ \\ \hline
$\phi _{gh^{k}}$ & $1$ & $-1$ \\ \hline
$\phi _{h^{k}}$ & $-1$ & $-1$ \\ \hline
\end{tabular}

Table 7.

Using Table 7 and (\ref{SLT}) it is easy to check that $\Gamma
_{h,g}=\langle \Gamma _{0}^{0}(2k,2),\gamma _{2}^{(k)}\rangle $ of genus
zero with $\gamma _{2}^{(k)}=\left( 
\begin{array}{cc}
a & b \\ 
ck & d
\end{array}
\right) $, $a=0$ mod $2$, $b,c,d\neq 0$ mod $2$; $\det \gamma _{2}^{(k)}=1$. 
$\gamma _{2}^{(k)}$ is of order $3$ in and normalises $\Gamma _{0}^{0}(2k,2)$%
. Then $\tilde{\Gamma}_{h,g}=2k|2-$.

Let us now examine the $p=3$ case. Since $g$, $h^{k}$ and $gh^{k}$ are all
Fricke and of order $3$ from Theorem \ref{theorem_h_p} IV.A we have $\phi
_{gh^{k}}(g)=\ \phi _{gh^{k}}(h^{k})=\omega _{3}^{2}$ and $\phi
_{h^{k}}(g)=1 $. Then using $\phi _{u}(h^{k})=(\phi _{u}(h))^{k}\,$ for $%
u=gh^{k},h^{k}$ we obtain the following residues for the singular cusps

\begin{tabular}{|l|l|l|l|}
\hline
Phase & $g$ & $h^{k}$ & $h$ \\ \hline
$\phi _{g}$ & $\omega _{3}$ & $1$ & $1$ \\ \hline
$\phi _{gh^{k}}$ & $\omega _{3}^{2}$ & $\omega _{3}^{2}$ & $
\begin{array}{c}
\omega_{3}^{2}\text{ when }k=1\ \mathrm{mod}\ 3 \\ 
\omega_{3}\text{ when }k=-1\ \mathrm{mod}\ 3 
\end{array}
$ \\ \hline
$\phi _{h^{k}}$ & $1$ & $\omega _{3}$ & $
\begin{array}{c}
\omega _{3}\text{ when }k=1\ \mathrm{mod}\ 3 \\ 
\omega _{3}^{2}\text{ when }k=-1\ \mathrm{mod}\ { }3
\end{array}
$ \\ \hline
\end{tabular}

Table 8.

Using Table 8 and (\ref{SLT}) it is easy to check that $Z\left[ 
\begin{array}{c}
h \\ 
g
\end{array}
\right] (\gamma _{3i}^{(k)}\tau )=q^{-1/3}+0+O(q^{1/3})$ for $i=1$, $2$ with 
$\gamma _{31}^{(k)}\equiv T\gamma _{3,k}^{\pm 1}$,$\ \gamma
_{32}^{(k)}\equiv T^{2}\gamma _{3,k}^{\pm 2}$, when $k=\pm 1\ {\rm mod}\ 3$
and with $\gamma _{3,k}\equiv ST^{-k}S$. $\gamma _{3i}^{(k)}$ is of order $2$
in and normalises $\Gamma _{0}^{0}(3k,3)$. Following the standard argument
we then find that $\Gamma _{h,g}=\langle \Gamma _{0}^{0}(3k,3),\gamma
_{31}^{(k)},\gamma _{32}^{(k)}\rangle $ of genus zero where the cusps $%
\{i\infty ,\frac{1}{k},\frac{2}{k},\frac{3}{k}\}\,$are identified. Then $%
\tilde{\Gamma}_{h,g}=3k|3-$ for $k=2,5$ are the only possible such modular
groups for $k$ prime.

With $p=5$ since $g$, $h^{k}$ and $gh^{k}$ are all Fricke and of order $5$
then from Theorem \ref{theorem_h_p} IV.A and in a similar way to that shown
above we obtain the following residues for the singular cusps

\begin{tabular}{|l|l|l|l|}
\hline
Phase & $g$ & $h^{k}$ & $h$ \\ \hline
$\phi _{g}$ & $\omega _{5}$ & $1$ & $1$ \\ \hline
$\phi _{gh^{k}}$ & $\omega _{5}^{3}$ & $\omega _{5}^{3}$ & $
\begin{array}{c}
\omega _{5}^{\mp 2}\text{ when }k=\pm 1\ {\rm mod}\ 5 \\ 
\omega _{5}^{\mp 1}\text{ when }k=\pm 2\ {\rm mod}\ 5
\end{array}
$ \\ \hline
$\phi _{h^{k}}$ & $1$ & $\omega _{5}$ & $
\begin{array}{c}
\omega _{5}^{\pm 1}\text{ when }k=\pm 1\ {\rm mod}\ 5 \\ 
\omega _{5}^{\mp 2}\text{ when }k=\pm 2\ {\rm mod}\ 5
\end{array}
$ \\ \hline
\end{tabular}

Table 9.

Using Table 9 and (\ref{SLT}) we can show that $Z\left[ 
\begin{array}{c}
h \\ 
g
\end{array}
\right] (\gamma _{5i}^{(k)}\tau )=q^{-1/5}+0+O(q^{1/5})$ for $i=1$, $2$ with 
$\gamma _{51}^{(k)}=T^{\pm 2}\gamma _{5,k}$ when $k=\pm 1\ {\rm mod}\ 5$, $%
\gamma _{51}^{(k)}=T^{\pm 1}\gamma _{5,k}$ when $k=\pm 2\ {\rm mod}\ 5$ $\,$%
and $\gamma _{52}^{(k)}=T^{\mp 2}\gamma _{5,k}^{-1}$ when $k=\pm 1\ {\rm mod}\ 5$, $\gamma _{52}^{(k)}=T^{\mp 1}\gamma _{5,k}^{-1}$ when $k=\pm 2\ {\rm mod}\ 5$ where $\gamma _{5,k}=ST^{-k}S$. $\gamma _{5i}^{(k)}$ is of order $3$
in and normalises $\Gamma _{0}^{0}(5k,5)$. Following the standard argument
we then find that $\Gamma _{h,g}=\langle \Gamma _{0}^{0}(5k,5),\gamma
_{51}^{(k)},\gamma _{52}^{(k)}\rangle $ of genus zero where the cusps $%
\{i\infty ,\frac{1}{k},\frac{2}{k},\frac{3}{k},\frac{4}{k},\frac{5}{k}\}\,$%
are identified. However no such genus zero modular group exists for $k$ $\,$%
prime.

\textbf{\ IX.C}\textit{\ }In this case $h$ is a $pk-$ element which is
possible only if $h=6-$ or $10-$. So there are four possibilities: (1) $p=3$%
, $k=2$ (2) $p=5$, $k=2$ (3) $p=2$, $k=3$ and (4) $p=2$, $k=5$. Following
detailed arguments similar to Theorem \ref{theorem_h_p^2} VII.B we find that
(1) and (2) lead to the contradictory property $\chi _{h^{p}}(\mathcal{P}%
_{\langle g,h\rangle })<0$.

Consider $p=2\,$with $k=3,5$. Since $h$ is non-Fricke then from Theorem \ref
{theorem_proj2} (ii) with $u=g$ Fricke and $v=h$ we have $\phi
_{gh}^{0}(g)\neq 1$ and hence $\phi _{gh}^{0}(g)=-1$. Consider $\gamma
_{2}^{(k)}=\left( 
\begin{array}{cc}
1 & 1-k \\ 
-1 & k
\end{array}
\right) \left( 
\begin{array}{cc}
k & 0 \\ 
0 & 1
\end{array}
\right) $ which normalises and is of order $2$ in $\Gamma _{0}^{0}(2k,2)$.
Then 
\[
Z\left[ 
\begin{array}{c}
h \\ 
g
\end{array}
\right] (\gamma _{2}^{(k)}\tau )=\phi
_{gh}^{0}((gh)^{k}g^{-1})q^{-1/2}+0+O(q^{1/2}). 
\]
However since $\phi _{gh}^{0}(gh)=\omega _{2k}$ from (\ref{FA}) and $\phi
_{gh}^{0}(g)=-1$ so we have $\phi _{gh}^{0}((gh)^{k}g^{-1})=1$. Hence $%
\Gamma _{h,g}=\langle \Gamma _{0}^{0}(2k,2),\gamma _{2}^{(k)}\rangle $ of
genus zero or $\tilde{\Gamma}_{h,g}=4k\frac{1}{2}+4k$ \cite{Q1}. Thus $%
\tilde{\Gamma}_{h,g}=12\frac{1}{2}+12$ in case (3) and $\tilde{\Gamma}%
_{h,g}=20\frac{1}{2}+20$ in case (4) both of genus zero.

\textbf{X.C} From Theorem \ref{theorem_impossible} I this class structure is
impossible.

\textbf{XI.C} Since $g,h^{k}$ $=p+$ and $gh^{k}=p-$ we find from Theorem \ref
{theorem_h_p} II.A that (1) $p=2$ (2) $p=3$ or (3) $p=5$. Let $\phi
_{gh}(g)=\omega _{3}^{\alpha }$ and $\phi _{gh}(h)=\omega _{pk}^{\beta }$.
From (\ref{FA}) $\phi _{gh}(gh)=\omega _{pk}$ and so $\alpha k+\beta =1\ 
{\rm mod}\ pk$. Since $\mathcal{V}_{\text{orb}}^{\langle g,h\rangle }\simeq 
\mathcal{V}^{\Lambda }$ using Theorem \ref{theorem_proj2} with $u=g$ Fricke
and $v=h$ we have $\phi _{gh}(g)\neq 1\,$so that $\alpha \neq 0\ {\rm mod}\
p$. Furthermore, with $u=h^{k}$ Fricke and $v=gh^{1-k}$ we find $\phi
_{gh}(h^{k})\neq 1$ and therefore $\beta \neq 0\ {\rm mod}\ p$.

(1) $p=2$. Since $\alpha \neq 0\ {\rm mod}\ 2$ we have $\alpha =1$ and $%
\beta =1-k=0\ {\rm mod}\ 2$ which is a contradiction to $\beta \neq 0\ 
{\rm mod}\ 2$ since $k\neq p=2$. Therefore this class structure is
impossible.

(2) $p=3$. Here $h$ is of type $3k+3$ so that $k=2$ or $7$ only. Since $%
h^{k}=3+$ we find from \ref{theorem_h_p} II.A that $\phi _{h^{k}}(g)=1$.
From (\ref{FA}) $\phi _{h^{k}}(h^{k})=\omega _{3}$ so that $\phi
_{h^{k}}(h)=\omega _{3}^{2}$ for $k=2\,$and $\phi _{h^{k}}(h)=\omega _{3}$
for $k=7$.

For $k=2$ the constraints $\alpha ,\beta \neq 0\ {\rm mod}\ 3$ imply that $%
\phi _{gh}(g)=\omega _{3}$, $\phi _{gh}(h)=\omega _{6}^{-1}$ so that $\phi
_{gh}(g^{-1}h^{-2})=1$. Consider $\gamma _{4}(\tau )=TST^{2}(2\tau )$. $%
\gamma _{4}\,$which normalises and is of order $4$ in $\Gamma _{0}^{0}(6,3)$
and acts on the cusps $\{i\infty ,1,2,3/2\}$ corresponding to the $%
\{g,gh,g^{2}h,h^{2}\}$ Fricke-twisted sectors in the following way: $i\infty
\rightarrow 2\rightarrow 3/2\rightarrow 1\rightarrow i\infty $. Then $%
Z\left[ 
\begin{array}{c}
h \\ 
g
\end{array}
\right] (\gamma _{4}\tau )=\ Z\left[ 
\begin{array}{c}
g^{-1}h^{-2} \\ 
gh
\end{array}
\right] (2\tau )=1.q^{-1/3}+0+O(q^{1/3})$ and $Z\left[ 
\begin{array}{c}
h \\ 
g
\end{array}
\right] (\gamma _{4}^{2}\tau )=Z\left[ 
\begin{array}{c}
gh^{3} \\ 
h^{2}
\end{array}
\right] (\tau )=\ 1.q^{-1/3}+0+O(q^{1/3})$ so that in the usual way we find
that $\Gamma _{h,g}=$ $\langle \Gamma _{0}^{0}(6,3),\gamma _{4}\rangle $ of
genus zero so that $\tilde{\Gamma}_{h,g}=\Gamma [18\symbol{126}a]\equiv
\langle \Gamma _{0}(18),T^{1/3}W_{18}T^{1/3}\rangle $ whose hauptmodul has $%
q $ expansion denoted by $18z\,$\cite{FMN} or $18\symbol{126}a\,$\cite{N2}.

For $k=7$ the constraints $\alpha ,\beta \neq 0\ {\rm mod}\ 3$ imply that $%
\phi _{gh}(g)=\omega _{3}^{2}$ and $\phi _{gh}(h)=\omega _{21}^{8}$ so that $%
\phi _{gh}(gh^{-7})=1$. Consider $\gamma _{1}=\frac{1}{\sqrt{7}}\left( 
\begin{array}{cc}
7 & -13 \\ 
-7 & 14
\end{array}
\right) $ and $\gamma _{2}=\frac{1}{\sqrt{7}}\left( 
\begin{array}{cc}
14 & -13 \\ 
-7 & 7
\end{array}
\right) $ both of which normalise and are of order $2$ in $\Gamma
_{0}^{0}(21,3)$. $\gamma _{1},\gamma _{2}$ act on the cusps $\{i\infty
,1,2,3/7\}$ corresponding to the $\{h,gh,g^{2}h,h^{7}\}$ Fricke-twisted
sectors in the following way: $\gamma _{1}$: $\{i\infty
,2\}\longleftrightarrow \{1,3/7\}$ and $\gamma _{2}$: $\{i\infty
,1\}\longleftrightarrow \{2,3/7\}$. Then in the usual way we find that $%
\Gamma _{h,g}=$ $\langle \Gamma _{0}^{0}(21,3),\gamma _{1},\gamma
_{2}\rangle $ of genus zero so that $\tilde{\Gamma}_{h,g}=$ $\Gamma [63%
\symbol{126}a]$ the modular group whose hauptmodul has $q$ expansion denoted
by $63\symbol{126}a\,$\cite{N2}.

(3) $p=5$. Similarly to Theorem \ref{theorem_h_p^2} VII.B we find the
contradictory property $\chi _{h^{5}}^{0}(\mathcal{P}_{\langle g,h\rangle
})<0$ so that this class structure is impossible.

\textbf{XII.C} From Theorem \ref{theorem_proj2} (i) $\mathcal{V}_{\text{orb}%
}^{\langle g,h\rangle }\simeq \mathcal{V}^{\Lambda }$. Then from (ii) with $%
u=g$ Fricke and $v=h$ it follows that $\phi _{gh}(g)=\omega _{p}^{s}$ for $%
s\neq 0\ {\rm mod}\ p$. Taking $u=g^{a}h^{k}$ Fricke for $a\neq 0\ {\rm mod}\ p$ and $v=h$ it also follows that $\phi _{gh}(g^{a}h^{k})\neq 1$. Let $%
\phi _{gh}(h^{k})=\omega _{p}^{r}$ then if$\ r\neq 0\ {\rm mod}\ p$ we can
find $a$ such that $as=-r\ {\rm mod}\ p\,$so that $\phi
_{gh}(g^{a}h^{k})=1\,$which is impossible. Hence $\phi _{gh}(h^{k})=1$.
Similarly, from Theorem \ref{theorem_proj2} (ii) with $u=g$ Fricke and $v=h$
it follows that $\phi _{gh^{p}}(g)\neq 1$ and taking $u=g^{a}h^{k}$ Fricke
for $a\neq 0\ {\rm mod}\ p$ and $v=h$ it follows that $\phi
_{gh^{p}}(g^{a}h^{k})\neq 1$. As above we then find that $\phi
_{gh^{p}}(h^{k})=1$. Finally, as in VI.C, we can also prove that $\phi
_{gh^{k}}(h)=1$.

Together these results can be used to show that $\Gamma _{h,g}=\langle
\Gamma _{0}^{0}(pk,p),\gamma _{p},\hat{W}_{k}\rangle $ of genus zero where $%
\gamma _{p}=ST^{k}S\,$ as defined in VI.C and $\hat{W}_{k}=\left( 
\begin{array}{cc}
ak & pb \\ 
ck & kd
\end{array}
\right) =\left( 
\begin{array}{cc}
a & pb \\ 
c & kd
\end{array}
\right) \left( 
\begin{array}{cc}
k & 0 \\ 
0 & 1
\end{array}
\right) $, $\det (\hat{W}_{k})=k$; $akd-bcp=1$. Therefore we find that $%
\tilde{\Gamma}_{h,g}=pk+k$.

\textbf{XIII.C} In this case all twisted sectors are Fricke and all cusps
are singular. This case is closely related to case VIII.C. Since $%
g,h^{k},gh^{k}=p+$ from Theorem \ref{theorem_h_p} IV.A we must have $p=2,3$
or $5$. Furthermore, for each such $p$ there is exactly one rational class
of order $pk$ in $G_{p+}$ not already identified in cases I.C to XII.C
above. We consider each $p$ in turn.

When $p=2$ by repeated use of Theorem \ref{theorem_proj1} (iii) and (\ref{FA}%
) we obtain the following residues which augment Table 7 above:\smallskip

\begin{tabular}{|l|l|l|l|}
\hline
Phase. & $g$ & $h$ & $h^{k}$ \\ \hline
$\phi _{g}$ & $-1$ & $1$ & $1$ \\ \hline
$\phi _{h}$ & $\pm 1$ & $\omega _{2k}$ & $-1$ \\ \hline
$\phi _{gh}$ & $\mp 1$ & $\mp \omega _{2k}$ & $\pm 1$ \\ \hline
$\phi _{gh^{2}}$ & $-1$ & $\mp (-\omega _{2k})^{(1-k)/2}$ & $\mp 1$ \\ \hline
$\phi _{h^{k}}$ & $-1$ & $-1$ & $-1$ \\ \hline
$\phi _{gh^{k}}$ & $1$ & $-1$ & $-1$ \\ \hline
\end{tabular}

Table 10. Phases for $p=2$.

Using Table 10 and following the usual argument we can then show that $%
\tilde{\Gamma}_{h,g}=2k|2+k$ or $2k|2+k^{\prime }$ where the GMF is either
fixed (upper signs in Table 10) or negated (lower signs in Table 10) by the
involution $\hat{w}_{k}:\tau \rightarrow -1/k\tau $.

When $p=3,5$ the phase residues are presented in Table 11 which augment
Tables 8 and 9 above.\smallskip

\begin{tabular}{|l|l|l|l|l|}
\hline
Phase & $g$ & $h^{k}$ & $h$ & Parameters \\ \hline
$\phi _{g}$ & $\omega _{p}$ & $1$ & $1$ &  \\ \hline
$\phi _{h}$ & $1$ & $\omega _{p}$ & $\omega _{pk}$ &  \\ \hline
$\phi _{gh}$ & $\omega _{p}^{\lambda }$ & $\omega _{p}^{1-k\lambda }$ & $%
\omega _{pk}^{1-k\lambda }$ & 
\begin{tabular}{l}
$\lambda =-k\ {\rm mod}\ p$, $p=3$ \\ 
$\lambda =-2k\ {\rm mod}\ 5$, $k=\pm 1\ {\rm mod}\ 5$, $p=5$ \\ 
$\lambda =2k\ {\rm mod}\ 5$, $k=\pm 2\ {\rm mod}\ 5$, $p=5$%
\end{tabular}
\\ \hline
$\phi _{gh^{p}}$ & $\omega _{p}^{\alpha }$ & $1$ & $\omega _{k}^{\beta ^{2}}$
& $\alpha k+\beta p=1$ \\ \hline
$\phi _{h^{k}}$ & $1$ & $\omega _{p}$ & $\omega _{p}^{\gamma }$ & $\gamma
k=1\ {\rm mod}\ p$ \\ \hline
$\phi _{gh^{k}}$ & $\omega _{p}^{(p+1)/2}$ & $\omega _{p}^{(p+1)/2}$ & $%
\omega _{p}^{\delta }$ & $\delta k=\frac{p+1}{2}\ {\rm mod}\ p$ \\ \hline
\end{tabular}

Table 11. Phases for $p=3,5$.

\smallskip These phase residues are determined by the use of Theorem \ref
{theorem_proj1} for $p=3,5$. For example for $p=5$ and $k=2\ {\rm mod}\ 5$
we have $\phi _{gh^{k+1}}(gh^{k})=\omega _{5}^{a}$ for some $a$ since $%
o(gh^{k})=5$. We will show that $a=0\ {\rm mod}\ 5$. Let $G=h^{k}$ and $%
H=g^{\alpha }h^{5\beta }$ where $\alpha k+5\beta =1$ so that $g=H^{k}$
\thinspace and $h=G^{\alpha }H^{5\beta }$. Then we have $\phi
_{gh^{k+1}}(gh^{k})=\phi _{H^{k}(G^{\alpha }H^{5\beta })^{k+1}}(H^{k}G)=\phi
_{G^{\alpha +1}H^{k+5\beta }}(GH^{k})\,$. Furthermore $\phi _{G}(H)=1$ using 
$\phi _{h^{k}}(g)=1\,$ from Theorem \ref{theorem_h_p} IV.A so that $H\in
G_{G}$. But $H$ belongs to the same unique rational class in $G_{G}=G_{5+}$
as $h$ does in $G_{g}$ and hence $Z\left[ 
\begin{array}{c}
h \\ 
g
\end{array}
\right] =Z\left[ 
\begin{array}{c}
H \\ 
G
\end{array}
\right] $. This implies that $\phi _{gh^{k+1}}(gh^{k})=\phi _{g^{\alpha
+1}h^{k+5\beta }}(gh^{k})\,$. But $k=2\ {\rm mod}\ 5$ implies that $\alpha
=3\ {\rm mod}\ 5$ and $k+5\beta =(1+3k)\ {\rm mod}\ 5k$. Hence $\omega
_{5}^{a}=\phi _{g^{-1}h^{1+3k}}(gh^{k})=$ $\phi _{gh^{k+1}}(g^{-1}h^{bk})\,$
where $b(1+3k)=(k+1)\ {\rm mod}\ 5k$ using (\ref{conjugation property}).
But $b=-1\ {\rm mod}\ 5$ so that $\omega _{5}^{a}=\phi
_{gh^{k+1}}(g^{-1}h^{-k})=\omega _{5}^{-a}$. Hence $a=$ $0\ {\rm mod}\ 5$
so that $\phi _{gh^{k+1}}(gh^{k})=1$. By use of (\ref{conjugation property})
and (\ref{FA}) we obtain $\phi _{gh}(g)=\omega _{5}^{-1}\,$and $\phi
_{gh}(h)=\omega _{5k}^{1+k}$. We can similarly determine the other entries
in Table 11.

For $p=3$ we may follow case VIII.C and use Table 11 to find in the usual
way that $\Gamma _{h,g}=\langle \Gamma _{0}^{0}(3k,3),\gamma
_{31}^{(k)},\gamma _{32}^{(k)},\hat{w}_{k}\rangle $ of genus zero i.e. $%
\tilde{\Gamma}_{h,g}=3k|3+$ for $k=1\ {\rm mod}\ 3$ which exists for $%
k=7,13 $ and $\tilde{\Gamma}_{h,g}=3k||3+$ for $k=-1\ {\rm mod}\ 3$ which
exists for $k=2,5$.

For $p=5$ case we may follow case VIII.C and use Table 11 to find in the
usual way that $\Gamma _{h,g}=\langle \Gamma _{0}^{0}(5k,5),\gamma
_{51}^{(k)},\gamma _{52}^{(k)},\hat{w}_{k}\rangle $ of genus zero i.e. $%
\tilde{\Gamma}_{h,g}=5k||5+$ for $k=\pm 2\ {\rm mod}\ 5$ which exists for $%
k=2,3,7$. $\tilde{\Gamma}_{h,g}$ does not exist for $k=\pm 1\ {\rm mod}\ 5$.

Table 12 summarises the results for cases I-XIII.C:

\begin{tabular}{|l|l|l|l|l|l|}
\hline
Cases & $\tilde{\Gamma}_{h,g}$ & $p=2$ & $p=3$ & $p=5$ & $p=7$ \\ \hline
I.C & $p^{2}k-$ & $12-$ & $18-$ & - & - \\ \hline
II.C & $p^{2}k+p^{2}k$ & $
\begin{tabular}{l}
$12+12,$ \\ 
$20+20$%
\end{tabular}
$ & $18+18$ & $50+50^{*}$ & - \\ \hline
III.C & $p^{2}k+p^{2}$ & \multicolumn{1}{|c|}{$
\begin{array}{cc}
12+4, & 20+4
\end{array}
$} & $18+9$ & - & - \\ \hline
IV.C & $
\begin{array}{c}
p^{2}k+k \\ 
4k+k^{\prime }\quad (p=2)
\end{array}
$ & $
\begin{array}{cc}
12+3 & 28+7 \\ 
12+3^{\prime } & 20+5^{\prime }
\end{array}
$ & $18+2$ & - & - \\ \hline
V.C & $p^{2}k+$ & $
\begin{tabular}{ll}
$12+,$ & $20+,$ \\ 
$28+,$ & $44+,$ \\ 
$92+$ & 
\end{tabular}
$ & $
\begin{array}{c}
18+ \\ 
45+
\end{array}
$ & $50+$ & - \\ \hline
VI.C & $pk-$ & $
\begin{array}{cc}
6-, & 10-
\end{array}
$ & $6-$ & $10-$ & - \\ \hline
VII.C & Impossible & - & - & - & - \\ \hline
VIII.C & $pk|p-$ & $
\begin{array}{cc}
6|2-, & 14|2-^{*}
\end{array}
$ & $
\begin{array}{c}
6|3- \\ 
15|3-^{*}
\end{array}
$ & - & - \\ \hline
IX.C & $4k\frac{1}{2}+4k$, $p=2$ & $
\begin{array}{c}
12\frac{1}{2}+12, \\ 
20\frac{1}{2}+20
\end{array}
$ & - & - & - \\ \hline
X.C & Impossible & - & - & - & - \\ \hline
XI.C & $\Gamma [p^{2}k\symbol{126}a]$ & - & $k=2,7$ & - & - \\ \hline
XII.C & $pk+k$ & $
\begin{tabular}{ll}
$6+3,$ & $10+5,$ \\ 
$14+7,$ & $22+11$ \\ 
$46+23,$ & 
\end{tabular}
$ & $
\begin{array}{c}
6+2 \\ 
15+5 \\ 
33+11
\end{array}
$ & $10+2$ & $21+3$ \\ \hline
XIII.C & $
\begin{array}{c}
2k|2+k \\ 
2k|2+k^{\prime } \\ 
pk|p+ \\ 
pk||p+
\end{array}
$ & 
\begin{tabular}{l}
$6|2+3$ \\ 
$6|2+3^{\prime }$ \\ 
$10|2+5^{\prime }$ \\ 
$14|2+7$ \\ 
$14|2+7^{\prime }$ \\ 
$22|2+11^{\prime }$ \\ 
$26|2+13^{\prime }$ \\ 
$34|2+17^{\prime }$ \\ 
$38|2+19^{\prime }$%
\end{tabular}
& $
\begin{array}{c}
6||3+ \\ 
15||3+ \\ 
21|3+ \\ 
39|3+
\end{array}
$ & $
\begin{array}{c}
10||5+ \\ 
15||5+ \\ 
35||5+
\end{array}
$ & - \\ \hline
\end{tabular}

Table 12. * indicates that no such GMF occurs.

\subsection{Case D: Anomalous Classes $h=pk|p+...$}

There are a number of GMFs for $g=p+$ for $p=2,3$ for which $o(h)=pk$ but $h$
is a member of an anomalous ${\bf M}$ class i.e. $h=pk|p+...$ Then the
hauptmodul property for the GMF can be demonstrated as follows. The Thompson
series $T_{h}(\tau )$ (\ref{Tg}) has the following property \cite{CN} :

\begin{equation}
\left[ T_{h}(\tau /p)\right] ^{p}=T_{h^{p}}(\tau )+\mathrm{const.}
\label{Po}
\end{equation}
and similarly for GMF functions we have \cite{T2}, \cite{T3}

\begin{equation}
\left[ Z\left[ 
\begin{array}{c}
h \\ 
g
\end{array}
\right] (\tau )\right] ^{p}=Z\left[ 
\begin{array}{c}
h^{p} \\ 
g
\end{array}
\right] (p\tau )+\mathrm{const.}  \label{Gpo}
\end{equation}
We then obtain the following examples of such GMFs:

\textbf{\ I.D} $p=2$ with $h=4|2-$ and $h^{2}=2-$. Then $Z\left[ 
\begin{array}{c}
h^{2} \\ 
g
\end{array}
\right] (2\tau )$ is a hauptmodul either for $\Gamma _{0}(4)$ when $%
gh^{2}=2- $ or for $\Gamma _{0}(2)$ when $gh^{2}=2+$ from Table 3. Therefore
from (\ref{Gpo}) $\Gamma _{h,g}=8|2-\,$or $4|2-$ where only the first
example arises in practice.

\textbf{\ II.D} $p=2$ with $h=4|2+$ and $h^{2}=2+$. Then $Z\left[ 
\begin{array}{c}
h^{2} \\ 
g
\end{array}
\right] (2\tau )$ is a hauptmodul for $\Gamma _{0}(4)+$ when $gh^{2}=2-$
from Table 3. Therefore $\Gamma _{h,g}=8|2+4$ or $8|2+4^{\prime }$. If $%
gh^{2}=2+$ then $Z\left[ 
\begin{array}{c}
h^{2} \\ 
g
\end{array}
\right] (2\tau )\,$is the hauptmodul for $2|2$ but no such genus zero
modular group commensurable with $SL(2,{\bf Z})$ exists fixing $\left(
Z\left[ 
\begin{array}{c}
h^{2} \\ 
g
\end{array}
\right] (2\tau )+\mathrm{const}\right) ^{1/2}$.

\textbf{\ III.D} $p=3\,$ with $h=3|3$ and $h^{3}=1$. Then $Z\left[ 
\begin{array}{c}
h^{3} \\ 
g
\end{array}
\right] (p\tau )$ is the hauptmodul for $\Gamma _{0}(3)+$ and therefore $%
\Gamma _{h,g}=9|3+$.

\textbf{\ IV.D} $p=3$ with $h=6|3-$ and $h^{3}=2-$ . Then $gh^{3}=6+3\,$
(since $(gh^{3})^{2}=3+$ and $(gh^{3})^{3}=2-$) so that $Z\left[ 
\begin{array}{c}
h^{3} \\ 
g
\end{array}
\right] (p\tau )$ is the hauptmodul for $\Gamma _{0}(6)+3$ and therefore $%
\Gamma _{h,g}=18|3+3$.

\textbf{\ V.D} $p=3$ with $h=21|3+$ and $h^{3}=7+$ . Then $gh^{3}=21+$ and $%
Z\left[ 
\begin{array}{c}
h^{3} \\ 
g
\end{array}
\right] (p\tau )$ is the hauptmodul for $\Gamma _{0}(21)+$ and therefore $%
\Gamma _{h,g}=$ $63|3+$.

\textbf{\ VI.D} $p=3$ with $h=39|3+$ and $h^{3}=13+$ . Then $gh^{3}=39+$ and 
$Z\left[ 
\begin{array}{c}
h^{3} \\ 
g
\end{array}
\right] (p\tau )$ is the hauptmodul for $\Gamma _{0}(39)+$ and therefore $%
\Gamma _{h,g}=117|3+$.

\subsection{\protect\smallskip Case E: Anomalous classes of type $4k|n+2,...$
associated with the Baby monster}

Consider $h=4k|n+2,...\in {\bf M}$ for $k=1$ or $k$ prime. There are eight
such classes in the Monster. Then $g=h^{2k}=2+$ and hence $h$ is of order $4k
$ in $C_{g}=\mathrm{2.B}$ but is of order $2k$ in $G_{g}=\mathrm{B}$. The
corresponding GMF is then 
\begin{eqnarray}
Z\left[ 
\begin{array}{c}
h \\ 
g
\end{array}
\right] (\tau ) &=&Z\left[ 
\begin{array}{c}
h \\ 
h^{2k}
\end{array}
\right] (\tau )  \nonumber \\
&=&\varepsilon (g,h;ST^{2k}S)Z\left[ 
\begin{array}{c}
h \\ 
1
\end{array}
\right] (ST^{2k}S\tau )
\end{eqnarray}
where $\varepsilon (g,h;ST^{2k}S)$ $\,$is a phase that must be present since 
$h$ is anomalous. This GMF is directly related to a standard Thompson series
and is therefore a hauptmodul. We have the following examples: $h=4|2+$, $%
8|4+$, $12|2+$, $12|2+2$, $20|2+$, $28|2+$, $52|2+$ and $68|2+$.

\section{Conclusions}

\smallskip We have shown how Generalised Moonshine can be understood within
an abelian orbifolding setting and have explicitly demonstrated the genus
zero property for rational Generalised Moonshine Functions (GMFs) arising in
a number of non-trivial cases. We have also discussed other aspects of
Generalised Moonshine such as properties of the head character expansion of
GMFs and constraints on the Monster classes for products of commuting
Monster elements. The orbifold methods developed in this paper can in
principle be extended to analyse all GMFs towards proving the genus zero
property in general. Examples of GMFs with irrational coefficients using
these methods appear in \cite{IT}.

\section{Acknowledgments}

We thank Simon Norton for providing us with invaluable information on
replicable series and for a number of useful discussions. We also thank
Geoffrey Mason and Rex Dark for useful discussions. We also acknowledge
funding from Enterprise Ireland under the Basic Research Grant Scheme.

\section{Appendix A. Modular Groups in Monstrous Moonshine}

In this appendix we the describe various modular groups associated with
Thompson series and Generalised Moonshine Functions (GMFs). Let $\Gamma ={SL}%
(2,{\bf Z})$ denote the full modular group. We define the following standard
subgroups of $\Gamma $: 
\begin{equation}
\Gamma _{0}(N)\equiv \left\{ \left( 
\begin{array}{cc}
a & b \\ 
c & d
\end{array}
\right) \in \Gamma ,\ c=0\text{ }{\rm mod}\ N\right\} ,  \tag{A.1}
\end{equation}

\begin{equation}
\Gamma _{0}^{0}(m,n)\equiv \left\{ \left( 
\begin{array}{cc}
a & b \\ 
c & d
\end{array}
\right) \in \Gamma \text{, }b=0\ {\rm mod}\ n\text{, }c=0\text{ }{\rm mod}%
\ m\right\} ,  \tag{A.2}
\end{equation}
\begin{equation}
\Gamma (m,n)\equiv \left\{ \left( 
\begin{array}{cc}
a & b \\ 
c & d
\end{array}
\right) \in \Gamma \text{, }\left( 
\begin{array}{cc}
a & b \\ 
c & d
\end{array}
\right) =\left( 
\begin{array}{cc}
1\ {\rm mod}\ n & 0\ {\rm mod}\ n\text{,} \\ 
0\text{ }{\rm mod}\ m & 1\text{ }{\rm mod}\ m
\end{array}
\right) \right\} \text{.}  \tag{A.3}
\end{equation}
We also use the notation $\Gamma _{0}^{0}(m)\equiv \Gamma _{0}^{0}(m,m)$ and 
$\Gamma (m)\equiv \Gamma (m,m)$. Note that $\Gamma _{0}(mn)\,$and $\Gamma
_{0}^{0}(m,n)$ are conjugate where $\Gamma _{0}(mn)=$ $\theta _{n}\Gamma
_{0}^{0}(m,n)\theta _{n}^{-1}$ with $\theta _{n}\equiv \left( 
\begin{array}{cc}
1 & 0 \\ 
0 & n
\end{array}
\right) $.

The normalizer $\mathcal{N}(\Gamma _{0}(N))=\left\{ \rho \in SL(2,{\bf R})%
\text{ }|\text{ }\rho \Gamma _{0}(N)\rho ^{-1}=\Gamma _{0}(N)\right\} $ is
also required to describe Monstrous Moonshine \cite{CN}. Let $h$ be an
integer where $h^{2}|N$, ($h^{2}$ divides $N$) and let $N=nh$. Then we
define the following sets of matrices.

$\Gamma _{0}(n|h)$: The group of matrices of the form 
\begin{equation}
\left( 
\begin{array}{cc}
a & \frac{b}{h} \\ 
cn & d
\end{array}
\right) \text{, }\det =1\text{,}  \tag{A.4}
\end{equation}
where $a$, $b$, $c$, $d$ $\in {\bf Z}$. For $h$ the largest divisor of $24$
for which $h^{2}|N$, $\Gamma _{0}(n|h)$ forms a subgroup of $\mathcal{N}%
(\Gamma _{0}(N))$. For $h=1$, $\Gamma _{0}(n|h)=\Gamma _{0}(n)$.

$W_{e}$: The set of matrices for a given positive integer $e$%
\begin{equation}
\left( 
\begin{array}{cc}
ae & b \\ 
cN & de
\end{array}
\right) \text{, }\det =e\text{, }e||N\text{,}  \tag{A.5}
\end{equation}
where $a$, $b$, $c$, $d$ $\in {\bf Z}$. $e||N$ denotes the property that $%
e|N $, and the $(e,N/e)=1$. The set $W_{e}$ forms a single coset of $\Gamma
_{0}(N)$ in $\mathcal{N}(\Gamma _{0}(N))$ with $W_{1}=\Gamma _{0}(N)$. It is
straightforward to show (up to scale factors) that

\begin{equation}
W_{e}^{2}=1\text{ mod }(\Gamma _{0}(N))\text{, }%
W_{e_{1}}W_{e_{2}}=W_{e_{2}}W_{e_{1}}=W_{e_{3}}\text{ mod }(\Gamma _{0}(N)),%
\text{ }  \tag{A.6}  \label{We}
\end{equation}
where $e_{3}=e_{1}e_{2}/(e_{1},e_{2})$. The coset $W_{e}$ is referred to as
an Atkin-Lehner (AL) involution for $\Gamma _{0}(N)$. The simplest example
is the Fricke involution $W_{N}$ with coset representative $\left( 
\begin{array}{cc}
0 & 1 \\ 
-N & 0
\end{array}
\right) $ which generates $\tau \rightarrow -1/N\tau $ and interchanges the
cusp points at $\tau =i\infty $ and $\tau =0$. For $e\neq n$ we can choose
the coset representative $\left( 
\begin{array}{cc}
e & b \\ 
N & de
\end{array}
\right) $ where $ed-bN/e=1$ which interchanges the cusp points at $\tau
=i\infty $ and $\tau =e/N$.

$w_{e}$: The set of matrices for a given positive integer $e$ of the form 
\begin{equation}
\left( 
\begin{array}{cc}
ae & \frac{b}{h} \\ 
cn & de
\end{array}
\right) \text{, }\det =e\text{, }e||\frac{n}{h}\text{,}  \tag{A.7}
\end{equation}
where $a$, $b$, $c$, $d$ $\in {\bf Z}$. The set $w_{e}$ is called an
Atkin-Lehner (AL) involution for $\Gamma _{0}(n|h)$. The properties (\ref{We}%
) are similarly obeyed by $w_{e}$ with $\Gamma _{0}(N)$ replaced by $\Gamma
_{0}(n|h)$.

$\mathcal{N}(\Gamma _{0}(N))$: The Normalizer of $\Gamma _{0}(N)$ in $SL(2,%
{\bf R})$ is constructed by adjoining to $\Gamma _{0}(n|h)$ all its AL
involutions $w_{e_{1}}$, $w_{e_{2}}$, ... where $h$ is the largest divisor
of 24 with $h^{2}|N$ and $N=nh$ \cite{CN}.

$\Gamma _{0}(n|h)+e_{1}$, $e_{2}$,... : This denotes the group obtained by
adjoining to $\Gamma _{0}(n|h)$ a particular subset of AL involutions $%
w_{e_{1}}$, $w_{e_{2}}$, ... and forms a subgroup of $\mathcal{N}(\Gamma
_{0}(N))$.

\smallskip $\Gamma _{0}(n|h)+$ : This denotes the group obtained by
adjoining to $\Gamma _{0}(n|h)$ all its AL involutions and forms a subgroup
of $\mathcal{N}(\Gamma _{0}(N))$.

Sometimes we abbreviate $\Gamma _{0}(n|h)+e_{1}$, $e_{2}$,... by $n|h+e_{1}$%
, $e_{2}$,... and $\Gamma _{0}(n|h)+$ by $n|h+$ (respectively $n+e_{1}$, $%
e_{2}$,... and $n+\,$for $h=1$).

$n||h+:$ This denotes the modular group containing $\Gamma _{0}(nh)$ having
index $h$ in $\theta _{h}(\Gamma _{0}(\frac{n}{h})+)\theta _{h}^{-1}$. This
conjugate contains the transformation $\tau \rightarrow \tau +m/h$ for all $%
m $ which takes the corresponding hauptmodul to $h$ distinct ones labelled
by the value of $m$, $0\leq m<h$ \cite{FMN}.

The following theorem \cite{Q1} gives us information about the cusp points
of $\Gamma _{0}(N)$:

\begin{thm}
\label{theorem_gamma0_cusps} For $N$ square-free, all $\Gamma _{0}(N)$
inequivalent cusps can be represented by the rational numbers $a/b$ such
that $b>0$, $b|N$, $(a,b)=1$ and $0<a<N/b$; two cusps $a/b$ and $a_{1}/b_{1}$
being $\Gamma _{0}(N)$ equivalent if and only if $b=b_{1}$ and $a=a_{1}$ mod 
$(b,N/b)$.
\end{thm}

\section{Appendix B. Identification of Generalised Moonshine Functions}

In this appendix we consider the explicit form of the head character
expansion of the GMF with $\phi _{g}(h)=1$

\begin{equation}
Z\left[ 
\begin{array}{c}
h \\ 
g
\end{array}
\right] (o(g)\tau )=\frac{1}{q}+0+\sum_{s=1}^\infty
a_{s}(h)q^{s}  \tag{B.1}  \label{HCF}
\end{equation}
for $g=p+$ for $p=2,3,5,7$ (see Table 1) and with $h\in C_{p+}\,$ of order $%
o(h)=kp$ for $k=1$ and $k$ prime. Here $a_{s}(h)\equiv $ $a_{g,s}(h)\,$of (%
\ref{Z(h,Fricke)}). From (\ref{chim}) we therefore expect that 
\begin{equation}
a_{s}(h)=a_{s-p}(h)+\tilde{\chi}_{s}(h)  \label{aformula}
\end{equation}
for some character $\tilde{\chi}_{s}$ of $C_{p+}$.

For $p=2,3$ we give the first 10 coefficients $a_{s}(h)$ of (\ref{HCF}) in
terms of the irreducible characters of $C_{p+}$ from the ATLAS \cite{CCNPW}.
The irreducible expansion for $p=5,7\,$appear in \cite{Q1}, \cite{Q2}$.$
Then using \cite{Q1}, \cite{Q2}, \cite{FMN} and \cite{N2} we may identify
the genus zero fixing group $\tilde{\Gamma}_{h,g}$ in each case considered.

When $g=2+$ then $C_{p+}=\mathrm{2.B}$ the double cover of the Baby Monster $%
\mathrm{B}$ then the first 10 head characters are:

\begin{eqnarray}
a_{1} &=&\chi _{1}+\chi _{2},\quad a_{2}=\chi _{185},\quad a_{3}=2\chi
_{1}+\chi _{2}+\chi _{3}+\chi _{4},\quad a_{4}=2\chi _{185}+\chi _{186}, 
\nonumber \\
a_{5} &=&3\chi _{1}+3\chi _{2}+2\chi _{3}+\chi _{4}+\chi _{6}+\chi
_{7},\quad a_{6}=4\chi _{185}+2\chi _{186}+\chi _{187,}  \nonumber \\
a_{7} &=&6\chi _{1}+5\chi _{2}+4\chi _{3}+3\chi _{4}+\chi _{5}+2\chi
_{6}+\chi _{7}+\chi _{8}+\chi _{9}+\chi _{10},  \nonumber \\
a_{8} &=&8\chi _{185}+4\chi _{186}+3\chi _{187}+\chi _{188},  \nonumber \\
a_{9} &=&8\chi _{1}+10\chi _{2}+7\chi _{3}+4\chi _{4}+2\chi _{5}+5\chi
_{6}+4\chi _{7}+2\chi _{8}+2\chi _{9}  \nonumber \\
&&+2\chi _{10}+\chi _{11}+\chi _{12}+\chi _{14}+\chi _{16}+\chi _{17}, 
\nonumber \\
a_{10} &=&14\chi _{185}+9\chi _{186}+7\chi _{187}+3\chi _{188}+\chi
_{189}+\chi _{192}.  \label{Babyheadchars}
\end{eqnarray}
Clearly the property (\ref{aformula}) is observed. Furthermore for $s$ odd, $%
a_{s}$ is a character for $\mathrm{B}$ whereas for $s$ even, $a_{s}$ is a
character for $\mathrm{2.B}$ for which $g$ is represented by $-1$ as
discussed in section 3.1 (iii).

When $g=3+$ then $C_{p+}=\mathrm{3.Fi}$, the triple cover of the Fischer
group $\mathrm{Fi}$, then the first 10 head characters are:

\begin{eqnarray}
a_{1} &=&\chi _{109},\quad a_{2}=\chi _{1}+\chi _{2},\quad a_{3}=\chi
_{109}+\chi _{110},\quad a_{4}=\chi _{109}+\chi _{110}+\chi _{112}, 
\nonumber \\
a_{5} &=&3\chi _{1}+2\chi _{2}+\chi _{3}+\chi _{8},\quad \quad a_{6}=2\chi
_{109}+2\chi _{110}+\chi _{112}+\chi _{113},  \nonumber \\
\quad a_{7} &=&3\chi _{109}+3\chi _{110}+\chi _{111}+\chi _{112}+\chi
_{113}+\chi _{114},  \nonumber \\
a_{8} &=&4\chi _{1}+5\chi _{2}+2\chi _{3}+\chi _{5}+2\chi _{8}+\chi
_{9}+\chi _{10}+\chi _{13},  \nonumber \\
a_{9} &=&5\chi _{109}+5\chi _{110}+\chi _{111}+2\chi _{112}+3\chi
_{113}+\chi _{114}+\chi _{116}+\chi _{119},  \nonumber \\
a_{10} &=&6\chi _{109}+6\chi _{110}+\chi _{111}+4\chi _{112}+4\chi
_{113}+2\chi _{114}+\chi _{115}+\chi _{116}  \nonumber \\
&&+\chi _{118}+\chi _{119}+\chi _{122}.
\end{eqnarray}
Clearly the property (\ref{aformula}) is again observed. Furthermore for $%
s=2\ {\rm }\ 3$, $a_{s}$ is a character for $\mathrm{Fi}$ otherwise $%
a_{s}$ is a character for $\mathrm{3.Fi}$, the triple cover as discussed in
section in Section 3.1 (iii).


\begin{thebibliography}{99}
\bibitem{DHVW}  L. Dixon, J. A. Harvey, C. C. Vafa and E. Witten, Nucl.
Phys. \textbf{B261} (1985) 678; Nucl. Phys. \textbf{B274} (1986) 285.

\bibitem{Go}  P. Goddard, Proceedings of the CIRM Luminy conference,
(World Scientific, Singapore, 1989).

\bibitem{DGM}  L. Dolan, P. Goddard and P. Montague, Comm. Math. Phys. 
\textbf{179} (1996) 61.

\bibitem{FLM2}  I. Frenkel, J. Lepowsky and A. Meurman, Vertex
Operator Algebras and the Monster, (Academic Press, New York, 1988).

\bibitem{Ka}  V. Kac, Vertex Operator Algebras for Beginners, University
Lecture Series, Vol. 10, (AMS, Boston, 1998).

\bibitem{MN}  A. Matsuo and K. Nagatomo, Math.Soc.Japan Memoirs, \textbf{%
4 }(1999) 1.

\bibitem{T1}  M. P. Tuite, Commun. Math. Phys. {\bf 146} (1992) 277.

\bibitem{T2}  M. P. Tuite, Commun. Math. Phys. {\bf 166} (1995) 495.

\bibitem{DLM1}  C. Dong, H. Li and G. Mason, Commun. Math. Phys. 
\textbf{214} (2000) 1.

\bibitem{FLM1}  I. Frenkel, J. Lepowsky and A. Meurman, Proc. Natl.
Acad. Sci. USA {\bf 81} (1984) 3256.

\bibitem{CN}  J.H. Conway and S.P. Norton, Bull. London Math. Soc. 
\textbf{11} (1979) 308.

\bibitem{B}  R. Borcherds, Invent. Math. \textbf{109} (1992) 405.

\bibitem{N1}  S. P. Norton, Proc. Symp. Pure Math. \textbf{47} (1987)
208.

\bibitem{T3}  M. P. Tuite, Contemp. Math. \textbf{193} (1996) 353.


\bibitem{Z}  Y. Zhu, J. Amer. Math. Soc. \textbf{9} (1996) 237.


\bibitem{Se}  J-P. Serre, A course in arithmetic, (Springer Verlag,
Berlin, 1978).

\bibitem{Va}  C. Vafa, Nucl. Phys. \textbf{B273} (1986) 592.

\bibitem{Hu}  Y-Z. Huang, Contemp. Math. \textbf{193} (1996) 123.

\bibitem{Gr}  R. Griess, Inv. Math. \textbf{68} (1982) 1.

\bibitem{DM}  C. Dong and G. Mason, U.C. Santa Cruz Preprint (1992).

\bibitem{DLM2}  C. Dong, H. Li and G. Mason, Contemp. Math. \textbf{193%
} (1996) 25.

\bibitem{Q1}  L. Queen, Some Relations Between Finite Groups, Lie Groups
and Modular Functions-Ph.D. Dissertation, (University of Cambridge,
Cambridge, 1980).

\bibitem{Q2}  L. Queen, Mathematics of Computation, \textbf{37} (1981)
547.

\bibitem{CCNPW}  J.H. Conway, R.T. Curtis, S.P. Norton, R.A. Parker
and R.A. Wilson, An Atlas of Finite Groups, (Clarendon Press, Oxford, 1985).


\bibitem{L}  W. Ledermann, Introduction to group characters, (Cambridge
University Press, Cambridge, 1989).

\bibitem{FMN}  D. Ford, J. MacKay and S.P. Norton, Commun. Algebra 
\textbf{22} (1994) 5175.

\bibitem[N2]{N2}  S. P. Norton, Private communication.

\bibitem[IT]{IT}  R. Ivanov and M.P. Tuite, math.QA/0202275.

\end{thebibliography}
\end{document}